\documentclass[11pt]{article}   

\usepackage{amsmath, amsthm, amsfonts, amssymb}
\usepackage{newtxtext,newtxmath} %<<<<< seems better than mathptmx
\usepackage{tikz}
\usetikzlibrary{arrows,shapes}
\definecolor{gold}{rgb}{1,.70,.0}   % <<<<<<<<<<<< color gold defined
\usepackage{float}                 % use in \begin{figure}[H] to enforce position

\let\OLDthebibliography\thebibliography
\renewcommand\thebibliography[1]{
  \OLDthebibliography{#1}
  \setlength{\parskip}{3pt}
  \setlength{\itemsep}{0pt plus 0.3ex}
}

\usepackage[hang,small,sf,bf]{caption}

\newcommand{\QED}{\ $\square$}

\def\nofigures{0}   %1 no figures, 0=figures. .8=some  <<<<<<<<<<<<<<<

\def\vvec#1#2{
{\renewcommand*{\arraystretch}{.7}
\begin{bmatrix}#1\\#2\end{bmatrix}}
}

\makeatletter
\newcommand*\bigcdot{\mathpalette\bigcdot@{.41}}
\newcommand*\bigcdot@[2]{\mathbin{\vcenter{\hbox{\scalebox{#2}{$\m@th#1\bullet$}}}}}
\makeatother
\makeatletter
\newcommand*\bbigcdot{\mathpalette\bigcdot@{.61}}
\newcommand*\bbigcdot@[2]{\mathbin{\vcenter{\hbox{\scalebox{#2}{$\m@th#1\bullet$}}}}}
\makeatother

\def\smalll{\footnotesize}

\begin{document}

\title{\bf Tessellations and  Descartes disk configurations}

\author{Jerzy Kocik                  %\thanks{support}
\\ \small Department of Mathematics
\\ \small Southern Illinois University, Carbondale, IL62901
\\ \small jkocik{@}siu.edu  }

%\date{\scriptsize (10 Dec 2018)}
\date{}

\maketitle

\vspace{-.3in}

\begin{abstract}
\noindent
An intriguing correspondence between certain finite planar tessellations 
and the Descartes circle arrangements is presented.  
This correspondence may be viewed as a visualization of the spinor structure underlying Descartes circles.
\\[5pt]
{\small 
\\[3pt]

} %end of s
\end{abstract}

%----------------------------------------------------------------------
\section{Dodecagonal tessellations}

%\rm
Ren\'e Descartes found a beautiful relation for four mutually tangent disks
(such arrangement is called {\bf Descartes configuration}): 
%$$
\begin{equation}
\label{eq:Descartes}
2\, (A^2 + B^2 + C^2 + D^2) \ = \  (A + B + C + D)^2
\end{equation}
%$$
where $A$, $B$, $C$, and $D$ are the curvatures of the disks , i.e., reciprocals of the radii.  
%Expressing (1) in terms of radii would lose its simplicity.  
He shared it with Elizabeth II, princes of Bohemia, in a form of a problem:
given three mutually tangent circles, find the fourth tangent to the three.
(For more on Descartes configuration and Apollonian disk packing, 
see \cite{GLM4,jk-m,LMW,Mel,N,IS}.)  
Since \eqref{eq:Descartes} is quadratic, there are two solutions; 
here is an example of such pair: 
$$
                   (\,2,\, 3, \, 6,\,  23\,)  \quad \hbox{and}\quad  (\,2,\,  3,\,  6,\,  -1\,)
$$
Finding other integral solutions requires certain parametrization  \cite{GLM1, jk-a}.
\\[-7pt]

It turns out that there is an intriguing duality between certain tessellations and Descartes arrangements of four circles: 
$$
\left.
\begin{array}{rr}
\hbox{\sf Areas in}\\
\hbox{\sf a tessellation}
\end{array}
\right\}
\quad
\longleftrightarrow
\quad
\left\{
\begin{array}{ll}
\hbox{\sf Curvatures}\\
\hbox{\sf in Descartes}\\
\hbox{\sf configuration}
\end{array}
\right.
$$
This duality is the main  topic of the present notes.
But there is a pleasant byproduct: one can easily construct an integral tessellation
and obtain  therefore effortlessly an integral Descartes disk configuration.
\\[-7pt]

Here is the construction.
Consider three vectors $\mathbf a$, $\mathbf b$, $\mathbf c$ in a two dimensional vector space $\mathbb R^2$, 
such that 
$$
\mathbf a+\mathbf b+\mathbf c = \mathbf 0
$$
%(They may also be viewed as three complex numbers.)  
Construct a tessellated dodecagon from this triple of vectors as described in Figure~\ref{fig:create}. 
It consists of 15 parallelograms, three of which are squares and the remaining twelve consists of pairs of congruent pieces. 
Here are the steps.

%==========================================
\def\maly{.37}
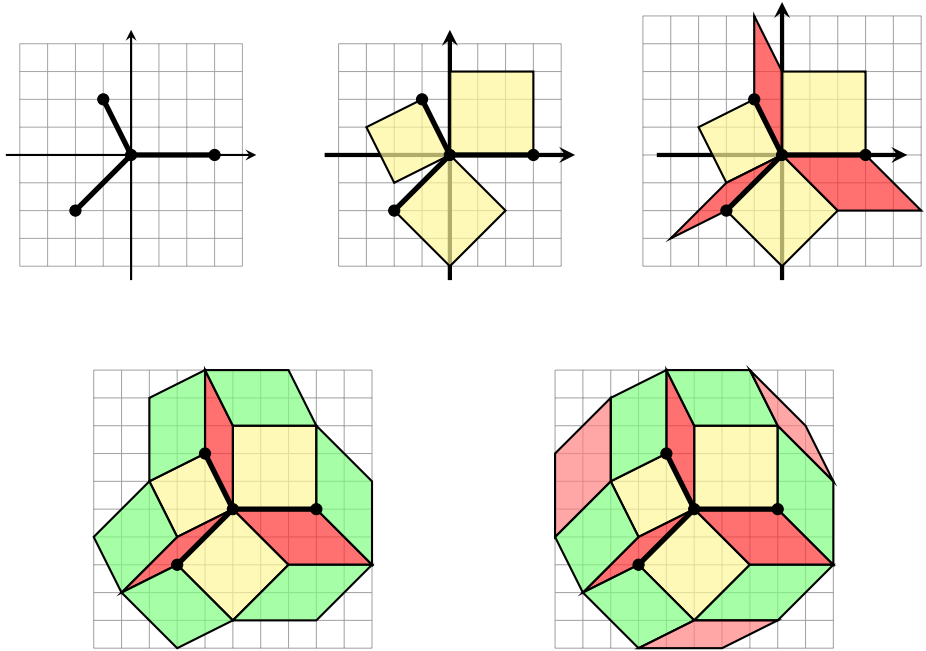
\begin{figure}[H]
\centering
%
%--------------------PANEL 1
\begin{tikzpicture}[scale=\maly, rotate=0, shift={(0,0cm)}]  %was 7
%\clip (-1.25,-2.1) rectangle (1.1,2.2);
\draw[very thin, color=gray!70] (-4,-4) grid (4,4);
\draw[-stealth, thick] (-4.5,0) -- (4.5,0);
\draw[-stealth, thick] (0, -4.5) -- (0, 4.5);

\def\Ax{3}; \def\Ay{0};
\def\Bx{-1}; \def\By{2};
\def\Cx{-2}; \def\Cy{-2};

\def\ax{0-\Ay}; \def\ay{\Ax}; \def\bx{0-\By}; \def\by{\Bx}; \def\cx{0-\Cy}; \def\cy{\Cx};

%lines and dots
\draw [fill=black] (0,0) circle  (.2);
\draw [fill=black]  (\Ax,\Ay) circle  (.2);
\draw [fill=black]  (\Bx,\By) circle  (.2);
\draw [fill=black]  (\Cx,\Cy) circle  (.2);
\draw [line width=2pt] (0,0) -- (\Ax,\Ay);
\draw [line width=2pt] (0,0) -- (\Bx,\By);
\draw [line width=2pt] (0,0) -- (\Cx,\Cy);

%\node at (3,1) [scale=1, color=black] {$\alpha$};
\draw[color=white, fill=white, fill opacity=\nofigures, draw opacity=\nofigures] (-4,-5) rectangle (4,5);     %<<<<<<<<<<<<<<<<<<<<<<<<<<<
\end{tikzpicture}
\qquad %=====================================
%
%%--------------------PANEL 2
\begin{tikzpicture}[scale=\maly, rotate=0, shift={(0,0cm)}]  %was 7

%\clip (-1.25,-2.1) rectangle (1.1,2.2);
\draw[thin, color=gray!70] (-4,-4) grid (4,4);
\draw[-stealth, ultra thick] (-4.5,0) -- (4.5,0);
\draw[-stealth, ultra thick] (0, -4.5) -- (0, 4.5);

\def\Ax{3}; \def\Ay{0};
\def\Bx{-1}; \def\By{2};
\def\Cx{-2}; \def\Cy{-2};
\def\ax{0-\Ay}; \def\ay{\Ax}; \def\bx{0-\By}; \def\by{\Bx}; \def\cx{0-\Cy}; \def\cy{\Cx};

%-----------------yellow
\draw [thick, fill=yellow!50, fill opacity=.7]  (0,0) -- (\Ax,\Ay) --  (\Ax+\ax, \Ay + \ay)  -- (\ax,\ay) -- cycle;
\draw [thick, fill=yellow!50, fill opacity=.7]  (0,0) -- (\Bx, \By) --  (\Bx+\bx, \By + \by)   -- (\bx,\by) -- cycle;
\draw [thick, fill=yellow!50, fill opacity=.7]  (0,0) -- (\Cx,\Cy) --  (\Cx+\cx,\Cy+\cy) -- (\cx,\cy) -- cycle;

%lines and dots
\draw [fill=black] (0,0) circle  (.2);
\draw [fill=black]  (\Ax,\Ay) circle  (.2);
\draw [fill=black]  (\Bx,\By) circle  (.2);
\draw [fill=black]  (\Cx,\Cy) circle  (.2);
\draw [line width=2pt] (0,0) -- (\Ax,\Ay);
\draw [line width=2pt] (0,0) -- (\Bx,\By);
\draw [line width=2pt] (0,0) -- (\Cx,\Cy);

%\node at (3,1) [scale=1, color=black] {$\alpha$};
%\draw [fill=gold!10] (.75,.75/2) circle  (.75/2);
\draw [color=white, fill=white, opacity=\nofigures, draw opacity=\nofigures] (-4,-5) rectangle (4,5);     %<<<<<<<<<<<<<<<<<<<<<<<<<<<
\end{tikzpicture}
\qquad   %=====================================================
\begin{tikzpicture}[scale=\maly, rotate=0, shift={(0,0cm)}]  %was 7

%\clip (-1.25,-2.1) rectangle (1.1,2.2);
\draw[thin, color=gray!70] (-5,-4) grid (5,5);
\draw[-stealth, ultra thick] (-4.5,0) -- (4.5,0);
\draw[-stealth, ultra thick] (0, -4.5) -- (0, 5.5);

\def\Ax{3}; \def\Ay{0};
\def\Bx{-1}; \def\By{2};
\def\Cx{-2}; \def\Cy{-2};
\def\ax{0-\Ay}; \def\ay{\Ax}; \def\bx{0-\By}; \def\by{\Bx}; \def\cx{0-\Cy}; \def\cy{\Cx};

%-----------------yellow
\draw [thick, fill=yellow!50, fill opacity=.7]  (0,0) -- (\Ax,\Ay) --  (\Ax+\ax, \Ay + \ay)  -- (\ax,\ay) -- cycle;
\draw [thick, fill=yellow!50, fill opacity=.7]  (0,0) -- (\Bx, \By) --  (\Bx+\bx, \By + \by)   -- (\bx,\by) -- cycle;
\draw [thick, fill=yellow!50, fill opacity=.7]  (0,0) -- (\Cx,\Cy) --  (\Cx+\cx,\Cy+\cy) -- (\cx,\cy) -- cycle;
%------------------red
\draw [thick, fill=red!80, fill opacity=.7]  (0,0) -- (\ax,\ay) -- (\ax+\Bx,\ay+\By) -- (\Bx,\By) -- cycle;
\draw [thick, fill=red!80, fill opacity=.7]  (0,0) -- (\bx,\by) -- (\bx+\Cx,\by+\Cy) -- (\Cx,\Cy) -- cycle;
\draw [thick, fill=red!80, fill opacity=.7]  (0,0) -- (\cx,\cy) -- (\cx+\Ax,\cy+\Ay) -- (\Ax,\Ay) -- cycle;

%lines and dots
\draw [fill=black] (0,0) circle  (.2);
\draw [fill=black]  (\Ax,\Ay) circle  (.2);
\draw [fill=black]  (\Bx,\By) circle  (.2);
\draw [fill=black]  (\Cx,\Cy) circle  (.2);
\draw [line width=2pt] (0,0) -- (\Ax,\Ay);
\draw [line width=2pt] (0,0) -- (\Bx,\By);
\draw [line width=2pt] (0,0) -- (\Cx,\Cy);

%\node at (3,1) [scale=1, color=black] {$\alpha$};
%\draw [fill=gold!10] (.75,.75/2) circle  (.75/2);
\draw [color=white, fill=white, opacity=\nofigures, draw opacity=\nofigures] (-5,-5) rectangle (5,5);     %<<<<<<<<<<<<<<<<<<<<<<<<<<<
\end{tikzpicture}

~\\~\\

\begin{tikzpicture}[scale=\maly, rotate=0, shift={(0,0cm)}]  %was 7
\draw[thin, color=gray!70] (-5,-5) grid (5,5);

%%\clip (-1.25,-2.1) rectangle (1.1,2.2);
%\draw (-\ile,0) -- (\ile,0);

\def\Ax{3}; \def\Ay{0};
\def\Bx{-1}; \def\By{2};
\def\Cx{-2}; \def\Cy{-2};
\def\ax{0-\Ay}; \def\ay{\Ax}; \def\bx{0-\By}; \def\by{\Bx}; \def\cx{0-\Cy}; \def\cy{\Cx};

%yellow
%-----------------yellow
\draw [thick, fill=yellow!50, fill opacity=.7]  (0,0) -- (\Ax,\Ay) --  (\Ax+\ax, \Ay + \ay)  -- (\ax,\ay) -- cycle;
\draw [thick, fill=yellow!50, fill opacity=.7]  (0,0) -- (\Bx, \By) --  (\Bx+\bx, \By + \by)   -- (\bx,\by) -- cycle;
\draw [thick, fill=yellow!50, fill opacity=.7]  (0,0) -- (\Cx,\Cy) --  (\Cx+\cx,\Cy+\cy) -- (\cx,\cy) -- cycle;
%------------------red
\draw [thick, fill=red!80, fill opacity=.7]  (0,0) -- (\ax,\ay) -- (\ax+\Bx,\ay+\By) -- (\Bx,\By) -- cycle;
\draw [thick, fill=red!80, fill opacity=.7]  (0,0) -- (\bx,\by) -- (\bx+\Cx,\by+\Cy) -- (\Cx,\Cy) -- cycle;
\draw [thick, fill=red!80, fill opacity=.7]  (0,0) -- (\cx,\cy) -- (\cx+\Ax,\cy+\Ay) -- (\Ax,\Ay) -- cycle;
%-----------------green
%aroun yellow 1, left, right
\draw [thick, fill=green!50, fill opacity=.7]  (\ax,\ay) -- (\Ax+\ax,\Ay+\ay) 
                                                                            -- (\Ax+\ax+\Bx,\Ay+\ay+\By) -- (\ax+\Bx,\ay+\By) -- cycle;
\draw [thick, fill=green!50, fill opacity=.7]  (\Ax,\Ay) -- (\Ax+\ax,\Ay+\ay) 
                                                                            -- (\Ax+\ax+\cx,\Ay+\ay+\cy) -- (\Ax+\cx,\Ay+\cy) -- cycle;                                                                   
%aroun yellow 2, left, right
\draw [thick, fill=green!50, fill opacity=.7]  (\bx,\by) -- (\Bx+\bx,\By+\by) 
                                                                            -- (\Bx+\bx+\Cx,\By+\by+\Cy) -- (\bx+\Cx,\by+\Cy) -- cycle;                                                                            
\draw [thick, fill=green!50, fill opacity=.7]  (\Bx,\By) -- (\Bx+\bx,\By+\by) 
                                                                            -- (\Bx+\bx+\ax,\By+\by+\ay) -- (\Bx+\ax,\By+\ay) -- cycle;                                                                   
%aroun yellow 3, left, right
\draw [thick, fill=green!50, fill opacity=.7]  (\cx,\cy) -- (\Cx+\cx,\Cy+\cy) 
                                                                            -- (\Cx+\cx+\Ax,\Cy+\cy+\Ay) -- (\cx+\Ax,\cy+\Ay) -- cycle;                                                                            
\draw [thick, fill=green!50, fill opacity=.7]  (\Cx,\Cy) -- (\Cx+\cx,\Cy+\cy) 
                                                                            -- (\Cx+\cx+\bx,\Cy+\cy+\by) -- (\Cx+\bx,\Cy+\by) -- cycle;                                                                                
                                                                          
%%--------------------lightred
%\draw [thick, fill=red!50, fill opacity=.7]   (\Ax+\ax,\Ay+\ay) -- (\Ax+\ax+\cx,\Ay+\ay+\cy)
%                                                        -- (\Ax+\ax+\Bx+\cx,\Ay+\ay+\By+\cy) -- (\Ax+\ax+\Bx,\Ay+\ay+\By)   -- cycle;                                                                            
%\draw [thick, fill=red!50, fill opacity=.7]   (\Bx+\bx,\By+\by) -- (\Bx+\bx+\ax,\By+\by+\ay)
%                                                        -- (\Bx+\bx+\Cx+\ax,\By+\by+\Cy+\ay) -- (\Bx+\bx+\Cx,\By+\by+\Cy)   -- cycle;                                                                            
%\draw [thick, fill=red!50, fill opacity=.7]   (\Cx+\cx,\cy+\cy) -- (\Cx+\cx+\bx,\Cy+\cy+\by)
%                                                        -- (\Cx+\cx+\Ax+\bx,\Cy+\cy+\Ay+\by) -- (\Cx+\cx+\Ax,\Cy+\cy+\Ay)   -- cycle;                                                                            
%lines and dots
\draw [fill=black] (0,0) circle  (.2);
\draw [fill=black]  (\Ax,\Ay) circle  (.2);
\draw [fill=black]  (\Bx,\By) circle  (.2);
\draw [fill=black]  (\Cx,\Cy) circle  (.2);
\draw [line width=2pt] (0,0) -- (\Ax,\Ay);
\draw [line width=2pt] (0,0) -- (\Bx,\By);
\draw [line width=2pt] (0,0) -- (\Cx,\Cy);

%\node at (3,1) [scale=1, color=black] {$\alpha$};
%\draw [fill=gold!10] (.75,.75/2) circle  (.75/2);
\draw [color=white, fill=white, opacity=\nofigures, draw opacity=\nofigures] (-5,-5) rectangle (5,5);     %<<<<<<<<<<<<<<<<<<<<<<<<<<<
\end{tikzpicture}
\qquad\qquad  \qquad%================================
\begin{tikzpicture}[scale=\maly, rotate=0, shift={(0,0cm)}]  %was 7

\draw[thin, color=gray!70] (-5,-5) grid (5,5);
%%\clip (-1.25,-2.1) rectangle (1.1,2.2);
%\draw (-\ile,0) -- (\ile,0);

\def\Ax{3}; \def\Ay{0};
\def\Bx{-1}; \def\By{2};
\def\Cx{-2}; \def\Cy{-2};
\def\ax{0-\Ay}; \def\ay{\Ax}; \def\bx{0-\By}; \def\by{\Bx}; \def\cx{0-\Cy}; \def\cy{\Cx};

%-----------------yellow
\draw [thick, fill=yellow!50, fill opacity=.7]  (0,0) -- (\Ax,\Ay) --  (\Ax+\ax, \Ay + \ay)  -- (\ax,\ay) -- cycle;
\draw [thick, fill=yellow!50, fill opacity=.7]  (0,0) -- (\Bx, \By) --  (\Bx+\bx, \By + \by)   -- (\bx,\by) -- cycle;
\draw [thick, fill=yellow!50, fill opacity=.7]  (0,0) -- (\Cx,\Cy) --  (\Cx+\cx,\Cy+\cy) -- (\cx,\cy) -- cycle;
%------------------red
\draw [thick, fill=red!80, fill opacity=.7]  (0,0) -- (\ax,\ay) -- (\ax+\Bx,\ay+\By) -- (\Bx,\By) -- cycle;
\draw [thick, fill=red!80, fill opacity=.7]  (0,0) -- (\bx,\by) -- (\bx+\Cx,\by+\Cy) -- (\Cx,\Cy) -- cycle;
\draw [thick, fill=red!80, fill opacity=.7]  (0,0) -- (\cx,\cy) -- (\cx+\Ax,\cy+\Ay) -- (\Ax,\Ay) -- cycle;
%-----------------green
%aroun yellow 1, left, right
\draw [thick, fill=green!50, fill opacity=.7]  (\ax,\ay) -- (\Ax+\ax,\Ay+\ay) 
                                                                            -- (\Ax+\ax+\Bx,\Ay+\ay+\By) -- (\ax+\Bx,\ay+\By) -- cycle;
\draw [thick, fill=green!50, fill opacity=.7]  (\Ax,\Ay) -- (\Ax+\ax,\Ay+\ay) 
                                                                            -- (\Ax+\ax+\cx,\Ay+\ay+\cy) -- (\Ax+\cx,\Ay+\cy) -- cycle;                                                                   
%aroun yellow 2, left, right
\draw [thick, fill=green!50, fill opacity=.7]  (\bx,\by) -- (\Bx+\bx,\By+\by) 
                                                                            -- (\Bx+\bx+\Cx,\By+\by+\Cy) -- (\bx+\Cx,\by+\Cy) -- cycle;                                                                            
\draw [thick, fill=green!50, fill opacity=.7]  (\Bx,\By) -- (\Bx+\bx,\By+\by) 
                                                                            -- (\Bx+\bx+\ax,\By+\by+\ay) -- (\Bx+\ax,\By+\ay) -- cycle;                                                                   
%aroun yellow 3, left, right
\draw [thick, fill=green!50, fill opacity=.7]  (\cx,\cy) -- (\Cx+\cx,\Cy+\cy) 
                                                                            -- (\Cx+\cx+\Ax,\Cy+\cy+\Ay) -- (\cx+\Ax,\cy+\Ay) -- cycle;                                                                            
\draw [thick, fill=green!50, fill opacity=.7]  (\Cx,\Cy) -- (\Cx+\cx,\Cy+\cy) 
                                                                            -- (\Cx+\cx+\bx,\Cy+\cy+\by) -- (\Cx+\bx,\Cy+\by) -- cycle;                                                                                
                                                                          
%--------------------lightred
\draw [thick, fill=red!50, fill opacity=.7]   (\Ax+\ax,\Ay+\ay) -- (\Ax+\ax+\cx,\Ay+\ay+\cy)
                                                        -- (\Ax+\ax+\Bx+\cx,\Ay+\ay+\By+\cy) -- (\Ax+\ax+\Bx,\Ay+\ay+\By)   -- cycle;                                                                            
\draw [thick, fill=red!50, fill opacity=.7]   (\Bx+\bx,\By+\by) -- (\Bx+\bx+\ax,\By+\by+\ay)
                                                        -- (\Bx+\bx+\Cx+\ax,\By+\by+\Cy+\ay) -- (\Bx+\bx+\Cx,\By+\by+\Cy)   -- cycle;                                                                            
\draw [thick, fill=red!50, fill opacity=.7]   (\Cx+\cx,\cy+\cy) -- (\Cx+\cx+\bx,\Cy+\cy+\by)
                                                        -- (\Cx+\cx+\Ax+\bx,\Cy+\cy+\Ay+\by) -- (\Cx+\cx+\Ax,\Cy+\cy+\Ay)   -- cycle;                                                                            
%lines and dots
\draw [fill=black] (0,0) circle  (.2);
\draw [fill=black]  (\Ax,\Ay) circle  (.2);
\draw [fill=black]  (\Bx,\By) circle  (.2);
\draw [fill=black]  (\Cx,\Cy) circle  (.2);
\draw [line width=2pt] (0,0) -- (\Ax,\Ay);
\draw [line width=2pt] (0,0) -- (\Bx,\By);
\draw [line width=2pt] (0,0) -- (\Cx,\Cy);

\draw [color=white, fill=white, opacity=\nofigures, draw opacity=\nofigures] (-5,-5) rectangle (5,5);     %<<<<<<<<<<<<<<<<<<<<<<<<<<<
\end{tikzpicture}

\caption{Creating a tassellation}
\label{fig:create}
\end{figure}

%Here is a more detailed description of Figure~\ref{fig:create}.
Start with three vectors with vanishing sum
%, $a+b+0$, 
 (first panel).  
Next, construct squares on these vectors (yellow squares in panel 2). 
Add parallelograms in spaces between the squares (red, panel 3).
Continue with adding the next layer of parallelograms as in Panel 4 (green).
Finally complete the construction with light-red parallelograms as in panel 5.
%The areas are easy to inspect.

\smallskip
 
 Now the unexpected happens:  interpret the areas of the three  red tiles at the center as curvatures of three disks,
 here 2, 3, and 6  (Figure~\ref{fig:two}). 
The curvature of the fourth disk, the solution to Descartes problem, 
can be read off from the above tessellation 
as the sum of the areas of three red and two green tiles, here 
 $2+3+2\times 6=23$.
 One may easily check that the Descartes formula holds:
 $$
 2(2^2+3^2+6^2+23^2)  =  (2+3+6+23)^2
 $$
 The curvature of the fourth disk may be read off from the picture in a few ways,
 for instance as the area of any of the three ``butterflies'' shown in Figure~\ref{fig:3butterflies}.

\begin{figure}[H]
\centering
\begin{tikzpicture}[scale=.45, rotate=0, shift={(0,0cm)}]  %was 7

\draw[thin, color=gray!70] (-5,-5) grid (5,5);
%%\clip (-1.25,-2.1) rectangle (1.1,2.2);
%\draw (-\ile,0) -- (\ile,0);

\def\Ax{3}; \def\Ay{0};
\def\Bx{-1}; \def\By{2};
\def\Cx{-2}; \def\Cy{-2};

\def\ax{0-\Ay}; \def\ay{\Ax}; \def\bx{0-\By}; \def\by{\Bx}; \def\cx{0-\Cy}; \def\cy{\Cx};

%-----------------yellow
\draw [thick, fill=yellow!50, fill opacity=.7]  (0,0) -- (\Ax,\Ay) --  (\Ax+\ax, \Ay + \ay)  -- (\ax,\ay) -- cycle;
\draw [thick, fill=yellow!50, fill opacity=.7]  (0,0) -- (\Bx, \By) --  (\Bx+\bx, \By + \by)   -- (\bx,\by) -- cycle;
\draw [thick, fill=yellow!50, fill opacity=.7]  (0,0) -- (\Cx,\Cy) --  (\Cx+\cx,\Cy+\cy) -- (\cx,\cy) -- cycle;
%------------------red
\draw [thick, fill=red!80, fill opacity=.7]  (0,0) -- (\ax,\ay) -- (\ax+\Bx,\ay+\By) -- (\Bx,\By) -- cycle;
\draw [thick, fill=red!80, fill opacity=.7]  (0,0) -- (\bx,\by) -- (\bx+\Cx,\by+\Cy) -- (\Cx,\Cy) -- cycle;
\draw [thick, fill=red!80, fill opacity=.7]  (0,0) -- (\cx,\cy) -- (\cx+\Ax,\cy+\Ay) -- (\Ax,\Ay) -- cycle;
%-----------------green
%aroun yellow 1, left, right
\draw [thick, fill=green!50, fill opacity=.7]  (\ax,\ay) -- (\Ax+\ax,\Ay+\ay) 
                                                                            -- (\Ax+\ax+\Bx,\Ay+\ay+\By) -- (\ax+\Bx,\ay+\By) -- cycle;
\draw [thick, fill=green!50, fill opacity=.7]  (\Ax,\Ay) -- (\Ax+\ax,\Ay+\ay) 
                                                                            -- (\Ax+\ax+\cx,\Ay+\ay+\cy) -- (\Ax+\cx,\Ay+\cy) -- cycle;                                                                   
%aroun yellow 2, left, right
\draw [thick, fill=green!50, fill opacity=.7]  (\bx,\by) -- (\Bx+\bx,\By+\by) 
                                                                            -- (\Bx+\bx+\Cx,\By+\by+\Cy) -- (\bx+\Cx,\by+\Cy) -- cycle;                                                                            
\draw [thick, fill=green!50, fill opacity=.7]  (\Bx,\By) -- (\Bx+\bx,\By+\by) 
                                                                            -- (\Bx+\bx+\ax,\By+\by+\ay) -- (\Bx+\ax,\By+\ay) -- cycle;                                                                   
%aroun yellow 3, left, right
\draw [thick, fill=green!50, fill opacity=.7]  (\cx,\cy) -- (\Cx+\cx,\Cy+\cy) 
                                                                            -- (\Cx+\cx+\Ax,\Cy+\cy+\Ay) -- (\cx+\Ax,\cy+\Ay) -- cycle;                                                                            
\draw [thick, fill=green!50, fill opacity=.7]  (\Cx,\Cy) -- (\Cx+\cx,\Cy+\cy) 
                                                                            -- (\Cx+\cx+\bx,\Cy+\cy+\by) -- (\Cx+\bx,\Cy+\by) -- cycle;                                                                                
                                                                          
%--------------------lightred
\draw [thick, fill=red!50, fill opacity=.7]   (\Ax+\ax,\Ay+\ay) -- (\Ax+\ax+\cx,\Ay+\ay+\cy)
                                                        -- (\Ax+\ax+\Bx+\cx,\Ay+\ay+\By+\cy) -- (\Ax+\ax+\Bx,\Ay+\ay+\By)   -- cycle;                                                                            
\draw [thick, fill=red!50, fill opacity=.7]   (\Bx+\bx,\By+\by) -- (\Bx+\bx+\ax,\By+\by+\ay)
                                                        -- (\Bx+\bx+\Cx+\ax,\By+\by+\Cy+\ay) -- (\Bx+\bx+\Cx,\By+\by+\Cy)   -- cycle;                                                                            
\draw [thick, fill=red!50, fill opacity=.7]   (\Cx+\cx,\cy+\cy) -- (\Cx+\cx+\bx,\Cy+\cy+\by)
                                                        -- (\Cx+\cx+\Ax+\bx,\Cy+\cy+\Ay+\by) -- (\Cx+\cx+\Ax,\Cy+\cy+\Ay)   -- cycle;                                                                            
%lines and dots
\draw [fill=black] (0,0) circle  (.2);
\draw [fill=black]  (\Ax,\Ay) circle  (.2);
\draw [fill=black]  (\Bx,\By) circle  (.2);
\draw [fill=black]  (\Cx,\Cy) circle  (.2);
\draw [line width=2pt] (0,0) -- (\Ax,\Ay);
\draw [line width=2pt] (0,0) -- (\Bx,\By);
\draw [line width=2pt] (0,0) -- (\Cx,\Cy);

\node at (1.5,1.5) [scale=1.7, color=black] {\sf 9};
\node at (0,-2) [scale=1.7, color=black] {\sf 8};
\node at (-1.5,.5) [scale=1.7, color=black] {\sf 5};
\node at (-.4,2.5) [scale=1.7, color=black] {\sf 3};
\node at (-1.92,-1.45) [scale=1.0, color=black] {\sf 2};
\node at (2.5,-1) [scale=1.7, color=black] {\sf 6};
\node at (4,.5) [scale=1.7, color=black] {\sf 6};
\node at (1,4) [scale=1.7, color=black] {\sf 6};
\node at (-2,3) [scale=1.7, color=black] {\sf 6};
\node at (-3.5,-1) [scale=1.7, color=black] {\sf 6};
\node at (-2,-3.5) [scale=1.7, color=black] {\sf 6};
\node at (2.5,-3) [scale=1.7, color=black] {\sf 6};

\node at (-4,1.5) [scale=1.7, color=black] {\sf 6};
\node at (.5,-4.53) [scale=1.5, color=black] {\sf 3};
\node at (3.5,3) [scale=1.2, color=black] {\sf 2};

%\draw [fill=gold!10] (.75,.75/2) circle  (.75/2);
\draw [color=white, fill=white, opacity=\nofigures, draw opacity=\nofigures] (-5,-5) rectangle (5,5);     %<<<<<<<<<<<<<<<<<<<<<<<<<<<
\end{tikzpicture}
\qquad
\begin{tikzpicture}[scale=1, rotate=0, shift={(0,0cm)}]
\node at (0,0) {~};
\node at (0,2) [scale =2] {$\Rightarrow$};
\end{tikzpicture}
\quad
\begin{tikzpicture}[scale=5.5, rotate=0]
\clip (-.35,1/5) rectangle (1,  1);
%\draw [fill=blue, opacity=.125, thick] (0,0) circle (1);
%-----------------------  horizontal two big -- x = 0 (Jas i Malgosia)
\draw [thick, fill=gold!10] (1/2,0) circle (1/2)
node  [scale=2]  at (1/2,1/4) {\sf 2}; %   {$\frac{\a,\b}{\c,\d}$}

\draw [thick, fill=gold!10] (0,2/3) circle (1/3)
node  [scale=2]  at (0,2/3) {\sf 3}; %   {$\frac{\a,\b}{\c,\d}$}

\draw [thick, fill=gold!10] (3/6, 4/6) circle (1/6)
node  [scale=2]  at (3/6,4/6) {\sf 6}; %   {$\frac{\a,\b}{\c,\d}$}

\draw [thick, fill=gold!10] (8/23, 12/23) circle (1/23)
node  [scale=.9]  at (8/23,12/23) {\sf 23}; %   {$\frac{\a,\b}{\c,\d}$}

%\draw [thick, dotted] (8/23, 12/23) circle (4/23);
\draw[color=white, fill=white, opacity=\nofigures, draw opacity=\nofigures] (-1,-1) rectangle (1,1);     %<<<<<<<<<<<<<<<<<<<<<<<<<<<
\end{tikzpicture}

%\end{changemargin}
\caption{Four disks in Descartes configuration from a tessellation}
\label{fig:two}
\end{figure}
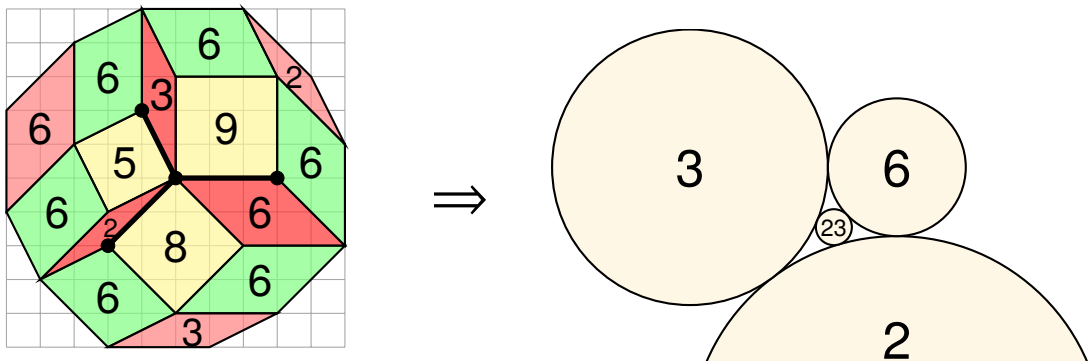

%============================FIGURE -- three butterflies =============
\begin{figure}
\centering 
\begin{tikzpicture}[scale=.25, rotate=0, shift={(0,0cm)}]  %was 7
\draw (-5,-5) grid (5,5);
%%\clip (-1.25,-2.1) rectangle (1.1,2.2);
%\draw (-\ile,0) -- (\ile,0);

\def\Ax{3}; \def\Ay{0};
\def\Bx{-1}; \def\By{2};
\def\Cx{-2}; \def\Cy{-2};
\def\ax{0-\Ay}; \def\ay{\Ax}; \def\bx{0-\By}; \def\by{\Bx}; \def\cx{0-\Cy}; \def\cy{\Cx};

%-----------------yellow
\draw [thick, fill=yellow!50, fill opacity=.7]  (0,0) -- (\Ax,\Ay) --  (\Ax+\ax, \Ay + \ay)  -- (\ax,\ay) -- cycle;
\draw [thick, fill=yellow!50, fill opacity=.7]  (0,0) -- (\Bx, \By) --  (\Bx+\bx, \By + \by)   -- (\bx,\by) -- cycle;
\draw [thick, fill=yellow!50, fill opacity=.7]  (0,0) -- (\Cx,\Cy) --  (\Cx+\cx,\Cy+\cy) -- (\cx,\cy) -- cycle;
%------------------red
\draw [thick, fill=red!80, fill opacity=.7]  (0,0) -- (\ax,\ay) -- (\ax+\Bx,\ay+\By) -- (\Bx,\By) -- cycle;
\draw [thick, fill=red!80, fill opacity=.7]  (0,0) -- (\bx,\by) -- (\bx+\Cx,\by+\Cy) -- (\Cx,\Cy) -- cycle;
\draw [thick, fill=red!80, fill opacity=.7]  (0,0) -- (\cx,\cy) -- (\cx+\Ax,\cy+\Ay) -- (\Ax,\Ay) -- cycle;
%-----------------green
%aroun yellow 1, left, right
\draw [thick, fill=green!50, fill opacity=.7]  (\ax,\ay) -- (\Ax+\ax,\Ay+\ay) 
                                                                            -- (\Ax+\ax+\Bx,\Ay+\ay+\By) -- (\ax+\Bx,\ay+\By) -- cycle;
\draw [thick, fill=green!50, fill opacity=.7]  (\Ax,\Ay) -- (\Ax+\ax,\Ay+\ay) 
                                                                            -- (\Ax+\ax+\cx,\Ay+\ay+\cy) -- (\Ax+\cx,\Ay+\cy) -- cycle;                                                                   
%aroun yellow 2, left, right
\draw [thick, fill=green!50, fill opacity=.7]  (\bx,\by) -- (\Bx+\bx,\By+\by) 
                                                                            -- (\Bx+\bx+\Cx,\By+\by+\Cy) -- (\bx+\Cx,\by+\Cy) -- cycle;                                                                            
\draw [thick, fill=green!50, fill opacity=.7]  (\Bx,\By) -- (\Bx+\bx,\By+\by) 
                                                                            -- (\Bx+\bx+\ax,\By+\by+\ay) -- (\Bx+\ax,\By+\ay) -- cycle;                                                                   
%aroun yellow 3, left, right
\draw [thick, fill=green!50, fill opacity=.7]  (\cx,\cy) -- (\Cx+\cx,\Cy+\cy) 
                                                                            -- (\Cx+\cx+\Ax,\Cy+\cy+\Ay) -- (\cx+\Ax,\cy+\Ay) -- cycle;                                                                            
\draw [thick, fill=green!50, fill opacity=.7]  (\Cx,\Cy) -- (\Cx+\cx,\Cy+\cy) 
                                                                            -- (\Cx+\cx+\bx,\Cy+\cy+\by) -- (\Cx+\bx,\Cy+\by) -- cycle;                                                                                
                                                                          
%--------------------lightred
\draw [thick, fill=red!50, fill opacity=.7]   (\Ax+\ax,\Ay+\ay) -- (\Ax+\ax+\cx,\Ay+\ay+\cy)
                                                        -- (\Ax+\ax+\Bx+\cx,\Ay+\ay+\By+\cy) -- (\Ax+\ax+\Bx,\Ay+\ay+\By)   -- cycle;                                                                            
\draw [thick, fill=red!50, fill opacity=.7]   (\Bx+\bx,\By+\by) -- (\Bx+\bx+\ax,\By+\by+\ay)
                                                        -- (\Bx+\bx+\Cx+\ax,\By+\by+\Cy+\ay) -- (\Bx+\bx+\Cx,\By+\by+\Cy)   -- cycle;                                                                            
\draw [thick, fill=red!50, fill opacity=.7]   (\Cx+\cx,\cy+\cy) -- (\Cx+\cx+\bx,\Cy+\cy+\by)
                                                        -- (\Cx+\cx+\Ax+\bx,\Cy+\cy+\Ay+\by) -- (\Cx+\cx+\Ax,\Cy+\cy+\Ay)   -- cycle;                                                                            

%lines and dots
\draw [fill=black] (0,0) circle  (.2);
\draw [fill=black]  (\Ax,\Ay) circle  (.2);
\draw [fill=black]  (\Bx,\By) circle  (.2);
\draw [fill=black]  (\Cx,\Cy) circle  (.2);
\draw [line width=2pt] (0,0) -- (\Ax,\Ay);
\draw [line width=2pt] (0,0) -- (\Bx,\By);
\draw [line width=2pt] (0,0) -- (\Cx,\Cy);

%opacity
\draw [white, fill=white, opacity=.749] (-5,-5) rectangle (5,5);

%===============shifted

%-----------------yellow
\draw [thick, fill=yellow!50, fill opacity=.9, xshift=1cm, yshift=1cm]  (0,0) -- (\Ax,\Ay) --  (\Ax+\ax, \Ay + \ay)  -- (\ax,\ay) -- cycle;
\draw [thick, fill=yellow!50, fill opacity=.9, xshift=-1.2cm, yshift=.6cm ]  (0,0) -- (\Bx, \By) --  (\Bx+\bx, \By + \by)   -- (\bx,\by) -- cycle;
\draw [thick, fill=yellow!50, fill opacity=.9, xshift=0cm, yshift=-1.4cm]  (0,0) -- (\Cx,\Cy) --  (\Cx+\cx,\Cy+\cy) -- (\cx,\cy) -- cycle;
%%------------------red
%\draw [thick, fill=red!80, fill opacity=.7]  (0,0) -- (\ax,\ay) -- (\ax+\Bx,\ay+\By) -- (\Bx,\By) -- cycle;
%\draw [thick, fill=red!80, fill opacity=.7]  (0,0) -- (\bx,\by) -- (\bx+\Cx,\by+\Cy) -- (\Cx,\Cy) -- cycle;
%\draw [thick, fill=red!80, fill opacity=.7]  (0,0) -- (\cx,\cy) -- (\cx+\Ax,\cy+\Ay) -- (\Ax,\Ay) -- cycle;
%-----------------green
%aroun yellow 1, left, right
\draw [thick, fill=green!50, fill opacity=.9, xshift=1cm, yshift=1cm]  (\ax,\ay) -- (\Ax+\ax,\Ay+\ay) 
                                                                            -- (\Ax+\ax+\Bx,\Ay+\ay+\By) -- (\ax+\Bx,\ay+\By) -- cycle;
\draw [thick, fill=green!50, fill opacity=.9, xshift=1cm, yshift=1cm]  (\Ax,\Ay) -- (\Ax+\ax,\Ay+\ay) 
                                                                            -- (\Ax+\ax+\cx,\Ay+\ay+\cy) -- (\Ax+\cx,\Ay+\cy) -- cycle;                                                                   
%%aroun yellow 2, left, right
\draw [thick, fill=green!50, fill opacity=.9, xshift=-1.2cm, yshift=.6cm]  (\bx,\by) -- (\Bx+\bx,\By+\by) 
                                                                            -- (\Bx+\bx+\Cx,\By+\by+\Cy) -- (\bx+\Cx,\by+\Cy) -- cycle;                                                                            
\draw [thick, fill=green!50, fill opacity=.9, xshift=-1.2cm, yshift=.6cm]  (\Bx,\By) -- (\Bx+\bx,\By+\by) 
                                                                            -- (\Bx+\bx+\ax,\By+\by+\ay) -- (\Bx+\ax,\By+\ay) -- cycle;                                                                   
%aroun yellow 3, left, right
\draw [thick, fill=green!50, fill opacity=.9, xshift=0cm, yshift=-1.4cm]  (\cx,\cy) -- (\Cx+\cx,\Cy+\cy) 
                                                                            -- (\Cx+\cx+\Ax,\Cy+\cy+\Ay) -- (\cx+\Ax,\cy+\Ay) -- cycle;                                                                            
\draw [thick, fill=green!50, fill opacity=.9, xshift=0cm, yshift=-1.4cm]  (\Cx,\Cy) -- (\Cx+\cx,\Cy+\cy) 
                                                                            -- (\Cx+\cx+\bx,\Cy+\cy+\by) -- (\Cx+\bx,\Cy+\by) -- cycle;                                                                                
%                                                                          
%--------------------lightred
\draw [thick, fill=red!50, fill opacity=.9, xshift=1cm, yshift=1cm]   (\Ax+\ax,\Ay+\ay) -- (\Ax+\ax+\cx,\Ay+\ay+\cy)
                                                        -- (\Ax+\ax+\Bx+\cx,\Ay+\ay+\By+\cy) -- (\Ax+\ax+\Bx,\Ay+\ay+\By)   -- cycle;                                                                            
\draw [thick, fill=red!50, fill opacity=.9, xshift=-1.2cm, yshift=.6cm]   (\Bx+\bx,\By+\by) -- (\Bx+\bx+\ax,\By+\by+\ay)
                                                        -- (\Bx+\bx+\Cx+\ax,\By+\by+\Cy+\ay) -- (\Bx+\bx+\Cx,\By+\by+\Cy)   -- cycle;                                                                            
\draw [thick, fill=red!50, fill opacity=.9, xshift=0cm, yshift=-1.4cm]   (\Cx+\cx,\cy+\cy) -- (\Cx+\cx+\bx,\Cy+\cy+\by)
                                                        -- (\Cx+\cx+\Ax+\bx,\Cy+\cy+\Ay+\by) -- (\Cx+\cx+\Ax,\Cy+\cy+\Ay)   -- cycle;                                                                            

\draw [ultra thick, xshift=1cm, yshift=1cm] (0,0) -- (\Ax,\Ay) -- (\Ax+\cx, \Ay+\cy) -- (\Ax+\ax+\cx, \Ay+\ay+\cy)
                                                               -- (\Ax+\ax+\Bx+\cx, \Ay+\ay+\By+\cy)  -- (\Ax+\ax+\Bx, \Ay+\ay+\By)
                                                               -- (\ax+\Bx,\ay+\By) --  (\ax,\ay)          -- cycle;

\draw [ultra thick, xshift=-1.2cm, yshift=.6cm] (0,0) -- (\Bx,\By) -- (\Bx+\ax, \By+\ay) -- (\Bx+\bx+\ax, \By+\by+\ay)
                                                               -- (\Bx+\bx+\Cx+\ax, \By+\by+\Cy+\ay)  -- (\Bx+\bx+\Cx, \By+\by+\Cy)
                                                               -- (\bx+\Cx,\by+\Cy) --  (\bx,\by)          -- cycle;

\draw [ultra thick, xshift=0cm, yshift=-1.4cm] (0,0) -- (\Cx,\Cy) -- (\Cx+\bx, \Cy+\by) -- (\Cx+\cx+\bx, \Cy+\cy+\by)
                                                               -- (\Cx+\cx+\Ax+\bx, \Cy+\cy+\Ay+\by)  -- (\Cx+\cx+\Ax, \Cy+\cy+\Ay)
                                                               -- (\cx+\Ax,\cy+\Ay) --  (\cx,\cy)          -- cycle;

%lines and dots
\draw [fill=black] (0,0) circle  (.2);
\draw [fill=black]  (\Ax,\Ay) circle  (.2);
\draw [fill=black]  (\Bx,\By) circle  (.2);
\draw [fill=black]  (\Cx,\Cy) circle  (.2);

\draw [line width=2pt] (0,0) -- (\Ax,\Ay);
\draw [line width=2pt] (0,0) -- (\Bx,\By);
\draw [line width=2pt] (0,0) -- (\Cx,\Cy);

%\node at (3,1) [scale=1, color=black] {$\alpha$};
%\draw [fill=gold!10] (.75,.75/2) circle  (.75/2);
\draw [draw =white, fill=white, fill opacity=\nofigures, draw opacity=\nofigures] (-7,-6.5) rectangle (7,6.1);     %<<<<<<<<<<<<<<<<<<<<<<<<<<<
\end{tikzpicture}
\caption{Three butterflies (the area of either equals to the curvature of the fourth disk)}
\label{fig:3butterflies}
\end{figure}
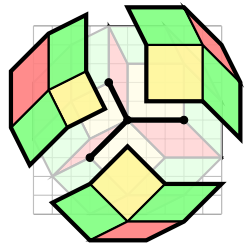

\noindent
{\bf Remark 1:}  The example above deals with integral vectors, which   
 makes the task of finding the areas easy.  In particular, one can use Pick’s Theorem \cite{Steinhaus}.

\smallskip
\smallskip
\noindent
{\bf Remark 2:}
 We witness here a non-intuitive duality. 
Since area is measured in $cm^2$ while curvatures in $cm^{–1}$,
the correspondence interchanges these units:    $cm^2\ \leftrightarrow \ cm^{–1}$.

\smallskip
\smallskip
\noindent
{\bf Observations:}  The following geometric facts are visible in the tessellation:  \\[-7pt]

1.	All green pieces have the same area (here equal to 6).

2.	Each shape among green pieces appears twice.

3.	The three exterior light-red pieces are congruent to the three red central pieces.

4.	Yellow squares have area equal to the sum of the adjacent red tiles.

5.  The sum of a square and  a red tile touching it at a vertex is the same for each square (here 11)

\begin{figure}[H]
\centering

\begin{tikzpicture}[scale=8.7, rotate=0]
\clip (0,1/3) rectangle (.7,  .91);
%\draw (0,0) circle (1);
%\draw [thick] (0,0) circle (1);
%\draw [fill=blue, opacity=.125, thick] (0,0) circle (1);
% A  1, 0, 2                  sums  12,  18,  34
% B  0, 2, 3
% C 3, 4, 6
%  D  8, 12 23 
\draw [thick] (1/2,0) circle (1/2)
node  [scale=2]  at (.62,.42) {\sf  2}; %   {$\frac{\a,\b}{\c,\d}$}

\draw [thick] (0,2/3) circle (1/3)
node  [scale=2]  at (1/50,2/3) {\sf 3}; %   {$\frac{\a,\b}{\c,\d}$}

\draw [thick] (3/6, 4/6) circle (1/6)
node  [scale=2]  at (3/6,4/6) {\sf 6}; %   {$\frac{\a,\b}{\c,\d}$}

% (A+B+C-D')/2   
\draw [thick] (8/23, 12/23) circle (1/23);
%node  [scale=.9]  at (8/23,12/23) {\sf 23}; %   {$\frac{\a,\b}{\c,\d}$}

% (A+B+C-D)/2   2, 3, 6
\draw [thick, dotted] (2/6, 3/6) circle (1/6);

% (A+B+D-C)/2     3, 5, 11
\draw [thick, dotted, fill=yellow!50, fill opacity=.70] (3/11, 5/11) circle (1/11);

% (A+C+D-B)/2     6, 7, 14
\draw [thick, dotted, fill=yellow!50, fill opacity=.70] (6/14, 7/14) circle (1/14);

% (B+C+D-A)/2     5, 9, 15
\draw [thick, dotted, fill=yellow!50, fill opacity=.70] (5/15, 9/15) circle (1/15);

\draw node  [scale=1.1]  at (.4,.37) {\color{red} \sf  6}; %   {$\frac{\a,\b}{\c,\d}$}
\draw node  [scale=1.1]  at (.24,.49) {\color{red} \sf 11}; %   {$\frac{\a,\b}{\c,\d}$}
\draw node  [scale=1.1]  at (.44,.54) {\color{red} \sf  14}; %   {$\frac{\a,\b}{\c,\d}$}
\draw node  [scale=1.1]  at (.30,.61) {\color{red} \sf  15}; %   {$\frac{\a,\b}{\c,\d}$}

\draw [color=white, fill=white, opacity=\nofigures, draw opacity=\nofigures] (-1,-1) rectangle (1,1);     %<<<<<<<<<<<<<<<<<<<<<<<<<<<
\end{tikzpicture}
\qquad\qquad %--------------------------------------------------------------------
\begin{tikzpicture}[scale=.45, rotate=0, shift={(0,0cm)}]  %was 7

\draw[thin, color=gray!70] (-5,-5) grid (5,5);
%%\clip (-1.25,-2.1) rectangle (1.1,2.2);
%\draw (-\ile,0) -- (\ile,0);

\def\Ax{3}; \def\Ay{0};
\def\Bx{-1}; \def\By{2};
\def\Cx{-2}; \def\Cy{-2};

\def\ax{0-\Ay}; \def\ay{\Ax}; \def\bx{0-\By}; \def\by{\Bx}; \def\cx{0-\Cy}; \def\cy{\Cx};

%-----------------yellow
\draw [thick, fill=yellow!50, fill opacity=.7]  (0,0) -- (\Ax,\Ay) --  (\Ax+\ax, \Ay + \ay)  -- (\ax,\ay) -- cycle;
\draw [thick, fill=yellow!50, fill opacity=.7]  (0,0) -- (\Bx, \By) --  (\Bx+\bx, \By + \by)   -- (\bx,\by) -- cycle;
\draw [thick, fill=yellow!50, fill opacity=.7]  (0,0) -- (\Cx,\Cy) --  (\Cx+\cx,\Cy+\cy) -- (\cx,\cy) -- cycle;
%------------------red
\draw [thick, fill=red!80, fill opacity=.7]  (0,0) -- (\ax,\ay) -- (\ax+\Bx,\ay+\By) -- (\Bx,\By) -- cycle;
\draw [thick, fill=red!80, fill opacity=.7]  (0,0) -- (\bx,\by) -- (\bx+\Cx,\by+\Cy) -- (\Cx,\Cy) -- cycle;
\draw [thick, fill=red!80, fill opacity=.7]  (0,0) -- (\cx,\cy) -- (\cx+\Ax,\cy+\Ay) -- (\Ax,\Ay) -- cycle;
%-----------------green
%aroun yellow 1, left, right
\draw [thick, fill=green!50, fill opacity=.7]  (\ax,\ay) -- (\Ax+\ax,\Ay+\ay) 
                                                                            -- (\Ax+\ax+\Bx,\Ay+\ay+\By) -- (\ax+\Bx,\ay+\By) -- cycle;
\draw [thick, fill=green!50, fill opacity=.7]  (\Ax,\Ay) -- (\Ax+\ax,\Ay+\ay) 
                                                                            -- (\Ax+\ax+\cx,\Ay+\ay+\cy) -- (\Ax+\cx,\Ay+\cy) -- cycle;                                                                   
%aroun yellow 2, left, right
\draw [thick, fill=green!50, fill opacity=.7]  (\bx,\by) -- (\Bx+\bx,\By+\by) 
                                                                            -- (\Bx+\bx+\Cx,\By+\by+\Cy) -- (\bx+\Cx,\by+\Cy) -- cycle;                                                                            
\draw [thick, fill=green!50, fill opacity=.7]  (\Bx,\By) -- (\Bx+\bx,\By+\by) 
                                                                            -- (\Bx+\bx+\ax,\By+\by+\ay) -- (\Bx+\ax,\By+\ay) -- cycle;                                                                   
%aroun yellow 3, left, right
\draw [thick, fill=green!50, fill opacity=.7]  (\cx,\cy) -- (\Cx+\cx,\Cy+\cy) 
                                                                            -- (\Cx+\cx+\Ax,\Cy+\cy+\Ay) -- (\cx+\Ax,\cy+\Ay) -- cycle;                                                                            
\draw [thick, fill=green!50, fill opacity=.7]  (\Cx,\Cy) -- (\Cx+\cx,\Cy+\cy) 
                                                                            -- (\Cx+\cx+\bx,\Cy+\cy+\by) -- (\Cx+\bx,\Cy+\by) -- cycle;                                                                                
                                                                          
%--------------------lightred
\draw [thick, fill=red!50, fill opacity=.7]   (\Ax+\ax,\Ay+\ay) -- (\Ax+\ax+\cx,\Ay+\ay+\cy)
                                                        -- (\Ax+\ax+\Bx+\cx,\Ay+\ay+\By+\cy) -- (\Ax+\ax+\Bx,\Ay+\ay+\By)   -- cycle;                                                                            
\draw [thick, fill=red!50, fill opacity=.7]   (\Bx+\bx,\By+\by) -- (\Bx+\bx+\ax,\By+\by+\ay)
                                                        -- (\Bx+\bx+\Cx+\ax,\By+\by+\Cy+\ay) -- (\Bx+\bx+\Cx,\By+\by+\Cy)   -- cycle;                                                                            
\draw [thick, fill=red!50, fill opacity=.7]   (\Cx+\cx,\cy+\cy) -- (\Cx+\cx+\bx,\Cy+\cy+\by)
                                                        -- (\Cx+\cx+\Ax+\bx,\Cy+\cy+\Ay+\by) -- (\Cx+\cx+\Ax,\Cy+\cy+\Ay)   -- cycle;

%lines and dots
\draw [fill=black] (0,0) circle  (.2);
\draw [fill=black]  (\Ax,\Ay) circle  (.2);
\draw [fill=black]  (\Bx,\By) circle  (.2);
\draw [fill=black]  (\Cx,\Cy) circle  (.2);
\draw [line width=2pt] (0,0) -- (\Ax,\Ay);
\draw [line width=2pt] (0,0) -- (\Bx,\By);
\draw [line width=2pt] (0,0) -- (\Cx,\Cy);

\node at (1.5,1.5) [scale=1.7, color=black] {\sf 9};
\node at (0,-2) [scale=1.7, color=black] {\sf 8};
\node at (-1.5,.5) [scale=1.7, color=black] {\sf 5};
\node at (-.4,2.5) [scale=1.7, color=black] {\sf 3};
\node at (-1.92,-1.45) [scale=1.1, color=black] {\sf 2};
\node at (2.5,-1) [scale=1.7, color=black] {\sf 6};
\node at (4,.5) [scale=1.7, color=black] {\sf 6};
\node at (1,4) [scale=1.7, color=black] {\sf 6};
\node at (-2,3) [scale=1.7, color=black] {\sf 6};
\node at (-3.5,-1) [scale=1.7, color=black] {\sf 6};
\node at (-2,-3.5) [scale=1.7, color=black] {\sf 6};
\node at (2.5,-3) [scale=1.7, color=black] {\sf 6};

\node at (-4,1.5) [scale=1.7, color=black] {\sf 6};
\node at (.5,-4.53) [scale=1.5, color=black] {\sf 3};
\node at (3.5,3) [scale=1.5, color=black] {\sf 2};

%opacity
\draw [white, fill=white, opacity=.89] (-5.1,-5.1) rectangle (5.1,5.1);

%-----------------yellow
\draw [thick, fill=yellow!50, fill opacity=.7]  (0,0) -- (\Ax,\Ay) --  (\Ax+\ax, \Ay + \ay)  -- (\ax,\ay) -- cycle;
\draw [thick, fill=yellow!50, fill opacity=.7]  (0,0) -- (\Bx, \By) --  (\Bx+\bx, \By + \by)   -- (\bx,\by) -- cycle;
\draw [thick, fill=yellow!50, fill opacity=.7]  (0,0) -- (\Cx,\Cy) --  (\Cx+\cx,\Cy+\cy) -- (\cx,\cy) -- cycle;

%-----------------green
%aroun yellow 1, left, right
\draw [thick, fill=green!50, fill opacity=.7]  (\ax,\ay) -- (\Ax+\ax,\Ay+\ay) 
                                                                            -- (\Ax+\ax+\Bx,\Ay+\ay+\By) -- (\ax+\Bx,\ay+\By) -- cycle;

\draw [thick, fill=green!50, fill opacity=.7]  (\bx,\by) -- (\Bx+\bx,\By+\by) 
                                                                            -- (\Bx+\bx+\Cx,\By+\by+\Cy) -- (\bx+\Cx,\by+\Cy) -- cycle;                                                                            

\draw [thick, fill=green!50, fill opacity=.7]  (\cx,\cy) -- (\Cx+\cx,\Cy+\cy) 
                                                                            -- (\Cx+\cx+\Ax,\Cy+\cy+\Ay) -- (\cx+\Ax,\cy+\Ay) -- cycle;                                                                            

\node at (1.5,1.5) [scale=1.7, color=black] {\sf 9};
\node at (0,-2) [scale=1.7, color=black] {\sf 8};
\node at (-1.5,.5) [scale=1.7, color=black] {\sf 5};

%\node at (4,.5) [scale=1.7, color=black] {\sf 6};
\node at (1,4) [scale=1.7, color=black] {\sf 6};
%\node at (-2,3) [scale=1.7, color=black] {\sf 6};
\node at (-3.5,-1) [scale=1.7, color=black] {\sf 6};
%\node at (-2,-3.5) [scale=1.7, color=black] {\sf 6};
\node at (2.5,-3) [scale=1.7, color=black] {\sf 6};

%\draw [fill=gold!10] (.75,.75/2) circle  (.75/2);
\draw [color=white, fill=white, opacity=\nofigures, draw opacity=\nofigures] (-5,-5) rectangle (5,5);     %<<<<<<<<<<<<<<<<<<<<<<<<<<<
\end{tikzpicture}

%\end{changemargin}
\caption{Mid-circles: \ 5+6=11, \ 8+6=14, \ 9+6=15.}
\label{fig:3circles}
\end{figure}
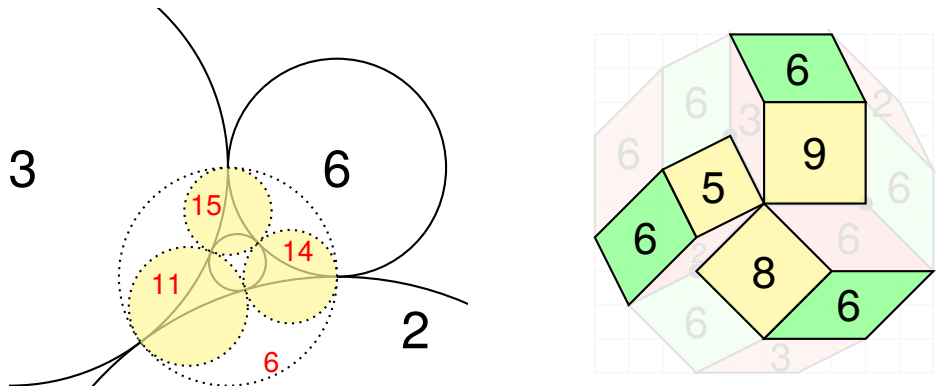

\smallskip
\smallskip
\noindent
There is more information on the Descartes configuration hidden in the tessellation.
A {\bf mid-circle} $(ABC)$ is defined as the unique circle that passes through the tangency points 
of a given 3 mutually tangent circles $A$, $B$, and $C$.\\[-17pt]
\begin{enumerate}
\item
The mid-circle of the three main disks (here 2, 3, 6) 
equals to the area the green tile  (here 6).\\[-19pt]
\item
The other mid-circles shown in Figure~\ref{fig:3circles} have curvatures equal
the area of a yellow tile plus the green, correspondingly.
(Here $5+6=11$, $8+6=14$, and $9+6=15$)
\end{enumerate}

The mid-circles relate to  symmetries of the original Descartes configuration
and its completion to the Apollonian Disk packing. 
For instance the circles 2, 3, 6 are invariant under inversions through circle~6.
The image of the smallest disk 23 is the external disk of curvature (-1).

There is, of course, the other solution of the Descartes problem.  
It may also be easily read off the tessellation,
and the corresponding mid-circles as well, namely
as a difference of the three red minus two green.
Or, equivalently, the sum of yellow and red that share point, minus two green.
(Here it is is $(-1)$, indicating an unbounded disk.)

\begin{figure}[H]
\centering

\begin{tikzpicture}[scale=3, rotate=0]
\clip (-.67,-.1) rectangle (1.6,  1.5);

\draw [fill=yellow!50, opacity=.7] (-1,-1) rectangle (2, 2);
\draw [white, fill=white] (1,1) circle (1);

\draw (0,0) circle (1);
%\draw [thick] (0,0) circle (1);
%\draw [fill=blue, opacity=.125, thick] (0,0) circle (1);
% A  1, 0, 2                  sums  12,  18,  34
% B  0, 2, 3
% C 3, 4, 6
%  D  8, 12 23 
\draw [thick] (1/2,0) circle (1/2)
node  [scale=2]  at (.40, -.02) {\sf  2}; %   {$\frac{\a,\b}{\c,\d}$}

\draw [thick] (0 ,2/3) circle (1/3)
node  [scale=2]  at (-.05,2/3) {\sf 3}; %   {$\frac{\a,\b}{\c,\d}$}

\draw [thick] (3/6, 4/6) circle (1/6)
node  [scale=2]  at (3/6-.05,4/6-.05) {\sf 6}; %   {$\frac{\a,\b}{\c,\d}$}

% (A+B+C-D')/2   
\draw [thick] (0, 0) circle (1);
\node  [scale=1.8]  at (-.5, 1.02) {\sf -1}; %   {$\frac{\a,\b}{\c,\d}$}

% (A+B+C-D)/2   2, 3, 6
\draw [thick, dotted] (1, 1) circle (1);

% (A+B+D-C)/2     3, 5, 11
\draw [thick, dotted, fill=yellow!50, fill opacity=.70] (1/3, 3/3) circle (1/3);

% (A+C+D-B)/2     6, 7, 14
\draw [black, thick, dotted, fill=yellow!50, fill opacity=.7] (2/2, 1/2) circle (1/2);

%% (B+C+D-A)/2     5, 9, 15
%\draw [thick, dotted, fill=yellow!50, fill opacity=.70] (5/15, 9/15) circle (1/15);

\draw node  [scale=1.1]  at (.2,1.1) {\color{red} \sf  3}; %   {$\frac{\a,\b}{\c,\d}$}
\draw node  [scale=1.1]  at (1.1,.49) {\color{red} \sf 2}; %   {$\frac{\a,\b}{\c,\d}$}
\draw node  [scale=1.1]  at (-.01,  1.42) {\color{red} \sf  -1}; %   {$\frac{\a,\b}{\c,\d}$}

\draw [color=white, fill=white, opacity=\nofigures, draw opacity=\nofigures] (-1,-1) rectangle (1.7,1.5);     %<<<<<<<<<<<<<<<<<<<<<<<<<<<
\end{tikzpicture}
\qquad\qquad %--------------------------------------------------------------------
\begin{tikzpicture}[scale=.45, rotate=0, shift={(0,0cm)}]  %was 7

\draw[thin, color=gray!70] (-5,-5) grid (5,5);
%%\clip (-1.25,-2.1) rectangle (1.1,2.2);
%\draw (-\ile,0) -- (\ile,0);

\def\Ax{3}; \def\Ay{0};
\def\Bx{-1}; \def\By{2};
\def\Cx{-2}; \def\Cy{-2};

\def\ax{0-\Ay}; \def\ay{\Ax}; \def\bx{0-\By}; \def\by{\Bx}; \def\cx{0-\Cy}; \def\cy{\Cx};

%-----------------yellow
\draw [thick, fill=yellow!50, fill opacity=.7]  (0,0) -- (\Ax,\Ay) --  (\Ax+\ax, \Ay + \ay)  -- (\ax,\ay) -- cycle;
\draw [thick, fill=yellow!50, fill opacity=.7]  (0,0) -- (\Bx, \By) --  (\Bx+\bx, \By + \by)   -- (\bx,\by) -- cycle;
\draw [thick, fill=yellow!50, fill opacity=.7]  (0,0) -- (\Cx,\Cy) --  (\Cx+\cx,\Cy+\cy) -- (\cx,\cy) -- cycle;
%------------------red
\draw [thick, fill=red!80, fill opacity=.7]  (0,0) -- (\ax,\ay) -- (\ax+\Bx,\ay+\By) -- (\Bx,\By) -- cycle;
\draw [thick, fill=red!80, fill opacity=.7]  (0,0) -- (\bx,\by) -- (\bx+\Cx,\by+\Cy) -- (\Cx,\Cy) -- cycle;
\draw [thick, fill=red!80, fill opacity=.7]  (0,0) -- (\cx,\cy) -- (\cx+\Ax,\cy+\Ay) -- (\Ax,\Ay) -- cycle;
%-----------------green
%aroun yellow 1, left, right
\draw [thick, fill=green!50, fill opacity=.7]  (\ax,\ay) -- (\Ax+\ax,\Ay+\ay) 
                                                                            -- (\Ax+\ax+\Bx,\Ay+\ay+\By) -- (\ax+\Bx,\ay+\By) -- cycle;
\draw [thick, fill=green!50, fill opacity=.7]  (\Ax,\Ay) -- (\Ax+\ax,\Ay+\ay) 
                                                                            -- (\Ax+\ax+\cx,\Ay+\ay+\cy) -- (\Ax+\cx,\Ay+\cy) -- cycle;                                                                   
%aroun yellow 2, left, right
\draw [thick, fill=green!50, fill opacity=.7]  (\bx,\by) -- (\Bx+\bx,\By+\by) 
                                                                            -- (\Bx+\bx+\Cx,\By+\by+\Cy) -- (\bx+\Cx,\by+\Cy) -- cycle;                                                                            
\draw [thick, fill=green!50, fill opacity=.7]  (\Bx,\By) -- (\Bx+\bx,\By+\by) 
                                                                            -- (\Bx+\bx+\ax,\By+\by+\ay) -- (\Bx+\ax,\By+\ay) -- cycle;                                                                   
%aroun yellow 3, left, right
\draw [thick, fill=green!50, fill opacity=.7]  (\cx,\cy) -- (\Cx+\cx,\Cy+\cy) 
                                                                            -- (\Cx+\cx+\Ax,\Cy+\cy+\Ay) -- (\cx+\Ax,\cy+\Ay) -- cycle;                                                                            
\draw [thick, fill=green!50, fill opacity=.7]  (\Cx,\Cy) -- (\Cx+\cx,\Cy+\cy) 
                                                                            -- (\Cx+\cx+\bx,\Cy+\cy+\by) -- (\Cx+\bx,\Cy+\by) -- cycle;                                                                                
                                                                          
%--------------------lightred
\draw [thick, fill=red!50, fill opacity=.7]   (\Ax+\ax,\Ay+\ay) -- (\Ax+\ax+\cx,\Ay+\ay+\cy)
                                                        -- (\Ax+\ax+\Bx+\cx,\Ay+\ay+\By+\cy) -- (\Ax+\ax+\Bx,\Ay+\ay+\By)   -- cycle;                                                                            
\draw [thick, fill=red!50, fill opacity=.7]   (\Bx+\bx,\By+\by) -- (\Bx+\bx+\ax,\By+\by+\ay)
                                                        -- (\Bx+\bx+\Cx+\ax,\By+\by+\Cy+\ay) -- (\Bx+\bx+\Cx,\By+\by+\Cy)   -- cycle;                                                                            
\draw [thick, fill=red!50, fill opacity=.7]   (\Cx+\cx,\cy+\cy) -- (\Cx+\cx+\bx,\Cy+\cy+\by)
                                                        -- (\Cx+\cx+\Ax+\bx,\Cy+\cy+\Ay+\by) -- (\Cx+\cx+\Ax,\Cy+\cy+\Ay)   -- cycle;

%lines and dots
\draw [fill=black] (0,0) circle  (.2);
\draw [fill=black]  (\Ax,\Ay) circle  (.2);
\draw [fill=black]  (\Bx,\By) circle  (.2);
\draw [fill=black]  (\Cx,\Cy) circle  (.2);
\draw [line width=2pt] (0,0) -- (\Ax,\Ay);
\draw [line width=2pt] (0,0) -- (\Bx,\By);
\draw [line width=2pt] (0,0) -- (\Cx,\Cy);

\node at (1.5,1.5) [scale=1.7, color=black] {\sf 9};
\node at (0,-2) [scale=1.7, color=black] {\sf 8};
\node at (-1.5,.5) [scale=1.7, color=black] {\sf 5};
\node at (-.4,2.5) [scale=1.7, color=black] {\sf 3};
\node at (-1.92,-1.45) [scale=1.1, color=black] {\sf 2};
\node at (2.5,-1) [scale=1.7, color=black] {\sf 6};
\node at (4,.5) [scale=1.7, color=black] {\sf 6};
\node at (1,4) [scale=1.7, color=black] {\sf 6};
\node at (-2,3) [scale=1.7, color=black] {\sf 6};
\node at (-3.5,-1) [scale=1.7, color=black] {\sf 6};
\node at (-2,-3.5) [scale=1.7, color=black] {\sf 6};
\node at (2.5,-3) [scale=1.7, color=black] {\sf 6};

\node at (-4,1.5) [scale=1.7, color=black] {\sf 6};
\node at (.5,-4.53) [scale=1.5, color=black] {\sf 3};
\node at (3.5,3) [scale=1.5, color=black] {\sf 2};

%opacity
\draw [white, fill=white, opacity=.89] (-5.1,-5.1) rectangle (5.1,5.1);

%-----------------yellow
\draw [thick, fill=yellow!50, fill opacity=.7]  (0,0) -- (\Ax,\Ay) --  (\Ax+\ax, \Ay + \ay)  -- (\ax,\ay) -- cycle;
\draw [thick, fill=yellow!50, fill opacity=.7]  (0,0) -- (\Bx, \By) --  (\Bx+\bx, \By + \by)   -- (\bx,\by) -- cycle;
\draw [thick, fill=yellow!50, fill opacity=.7]  (0,0) -- (\Cx,\Cy) --  (\Cx+\cx,\Cy+\cy) -- (\cx,\cy) -- cycle;

%-----------------green
\draw [thick, fill=green!50, fill opacity=.7]  (\ax,\ay) -- (\Ax+\ax,\Ay+\ay) 
                                                                            -- (\Ax+\ax+\Bx,\Ay+\ay+\By) -- (\ax+\Bx,\ay+\By) -- cycle;
\draw [thick, fill=green!50, fill opacity=.7]  (\bx,\by) -- (\Bx+\bx,\By+\by) 
                                                                            -- (\Bx+\bx+\Cx,\By+\by+\Cy) -- (\bx+\Cx,\by+\Cy) -- cycle;                                                                            
\draw [thick, fill=green!50, fill opacity=.7]  (\cx,\cy) -- (\Cx+\cx,\Cy+\cy) 
                                                                            -- (\Cx+\cx+\Ax,\Cy+\cy+\Ay) -- (\cx+\Ax,\cy+\Ay) -- cycle;                                                                            

\node at (1.5,1.5) [scale=1.7, color=black] {\sf 9};
\node at (0,-2) [scale=1.7, color=black] {\sf 8};
\node at (-1.5,.5) [scale=1.7, color=black] {\sf 5};

\node at (1,4) [scale=1.7, color=black] {\sf 6};
\node at (-3.5,-1) [scale=1.7, color=black] {\sf 6};
\node at (2.5,-3) [scale=1.7, color=black] {\sf 6};

\draw [color=white, fill=white, opacity=\nofigures, draw opacity=\nofigures] (-5,-5) rectangle (5,5);     %<<<<<<<<<<<<<<<<<<<<<<<<<<<
\end{tikzpicture}

%\end{changemargin}
\caption{The fourth circle $D'=-1$ and the mid-circles: \ $5-6=-1$, \ $8-6=2$, \ $9-6=3$  (numbers in red font).}
\label{fig:3circles2}
\end{figure}

The next two figures are just a few more examples of tiling, drawn with the help of Cinderella \cite{cinderella}. 
In general, any starting vectors are admissible,
There is a slight problem that for some choices of the initial vectors,
some tiles will have negative areas and therefore will overlap with other tiles.

%===============================================
\begin{figure}[H]
\centering
% trim: left bottom right top
\includegraphics[scale=.42, trim={2cm 3.7cm 0 1cm},clip] {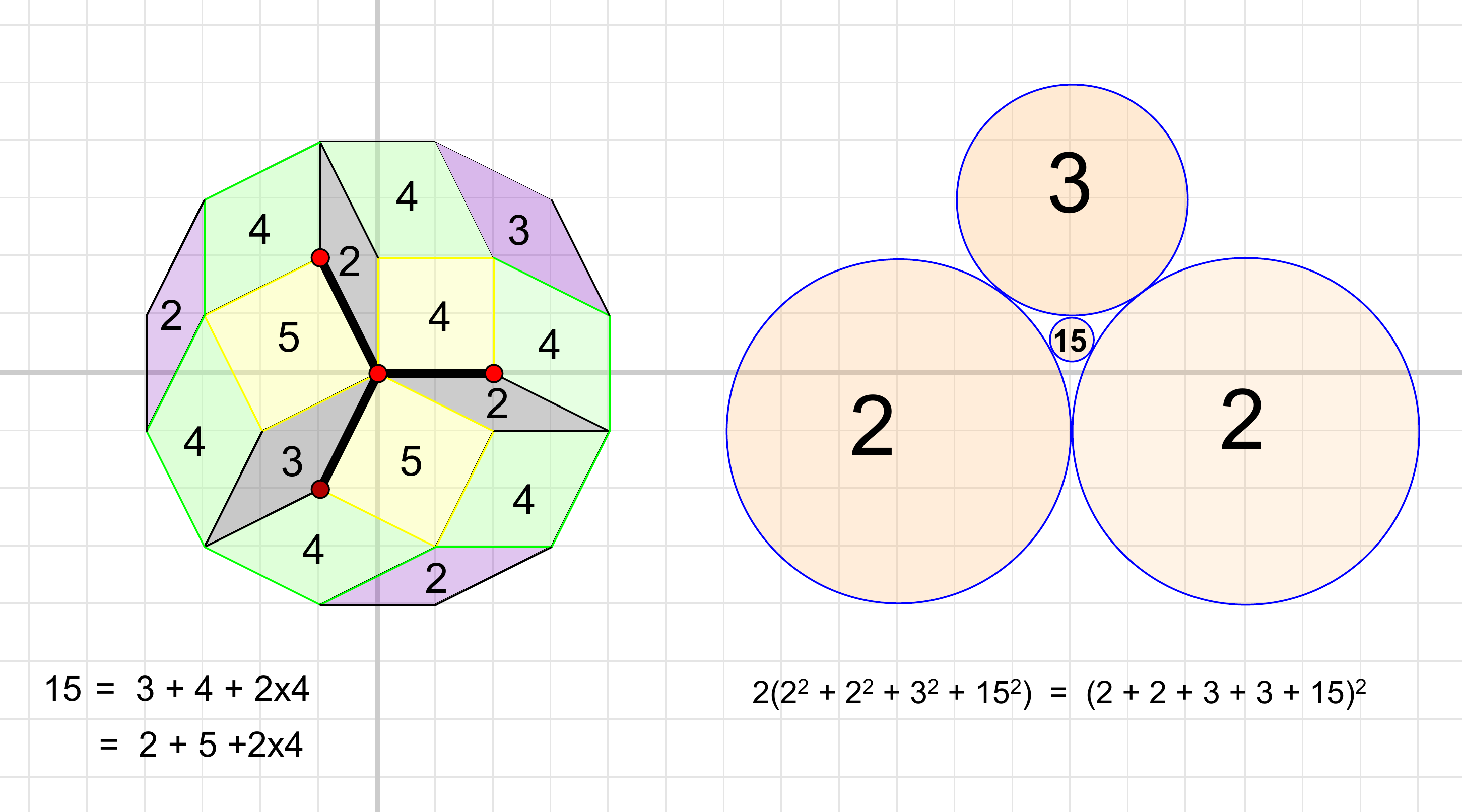}
\caption{Another example of tiling and the corresponding Descartes configuration: 
 $(2+3+3+15)^2$ \\ $= 2\, (\, 2^2+2^2+3^2+15^2\,) $.}
     %<<<<<<<<<<<<<<<<<<<<<<<<<<<
\label{fig:MainExample} 
\end{figure}

%===============================================
\begin{figure}[H]
\centering
\includegraphics[scale=.29] {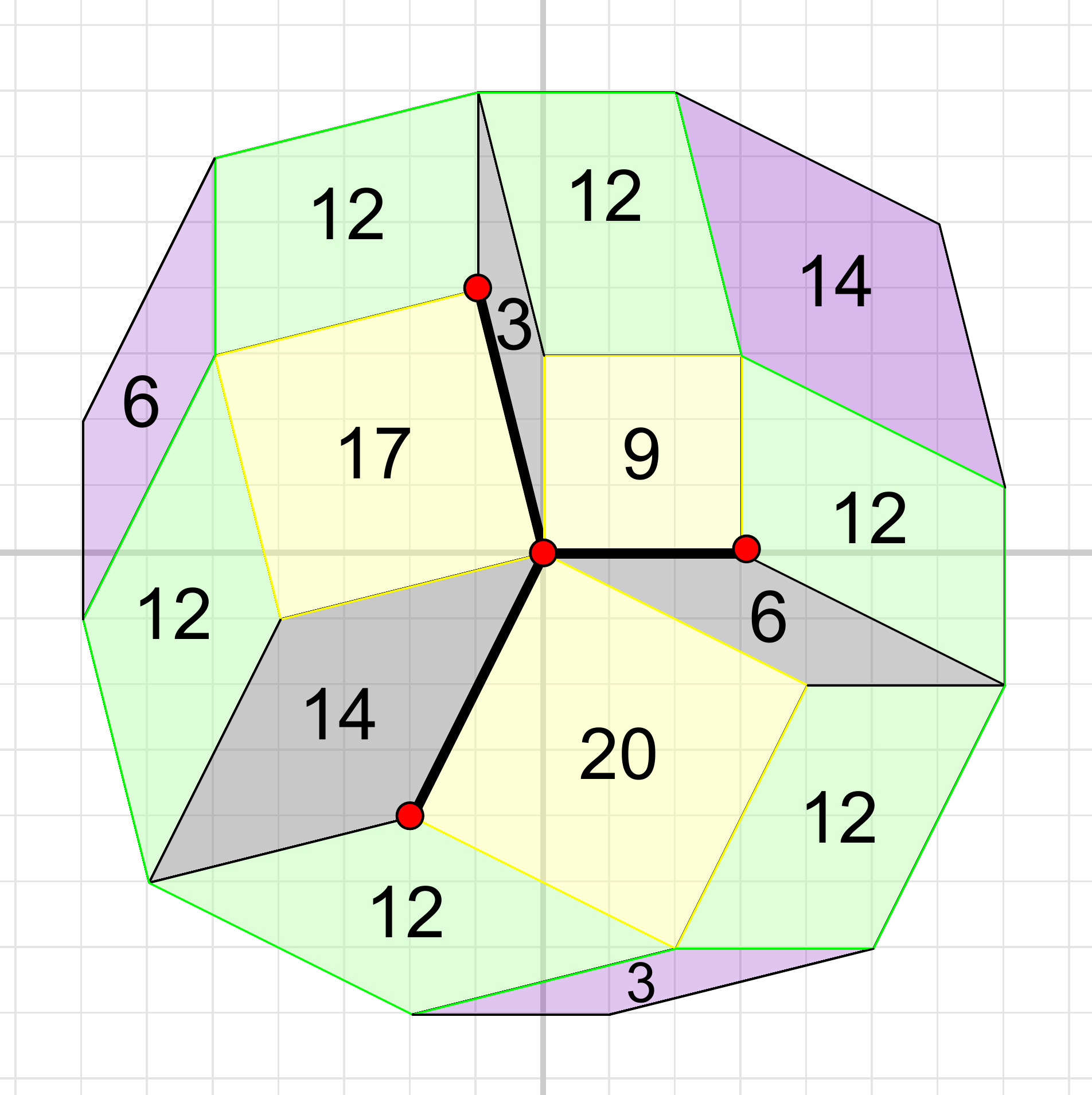}\quad 
\includegraphics[scale=.25] {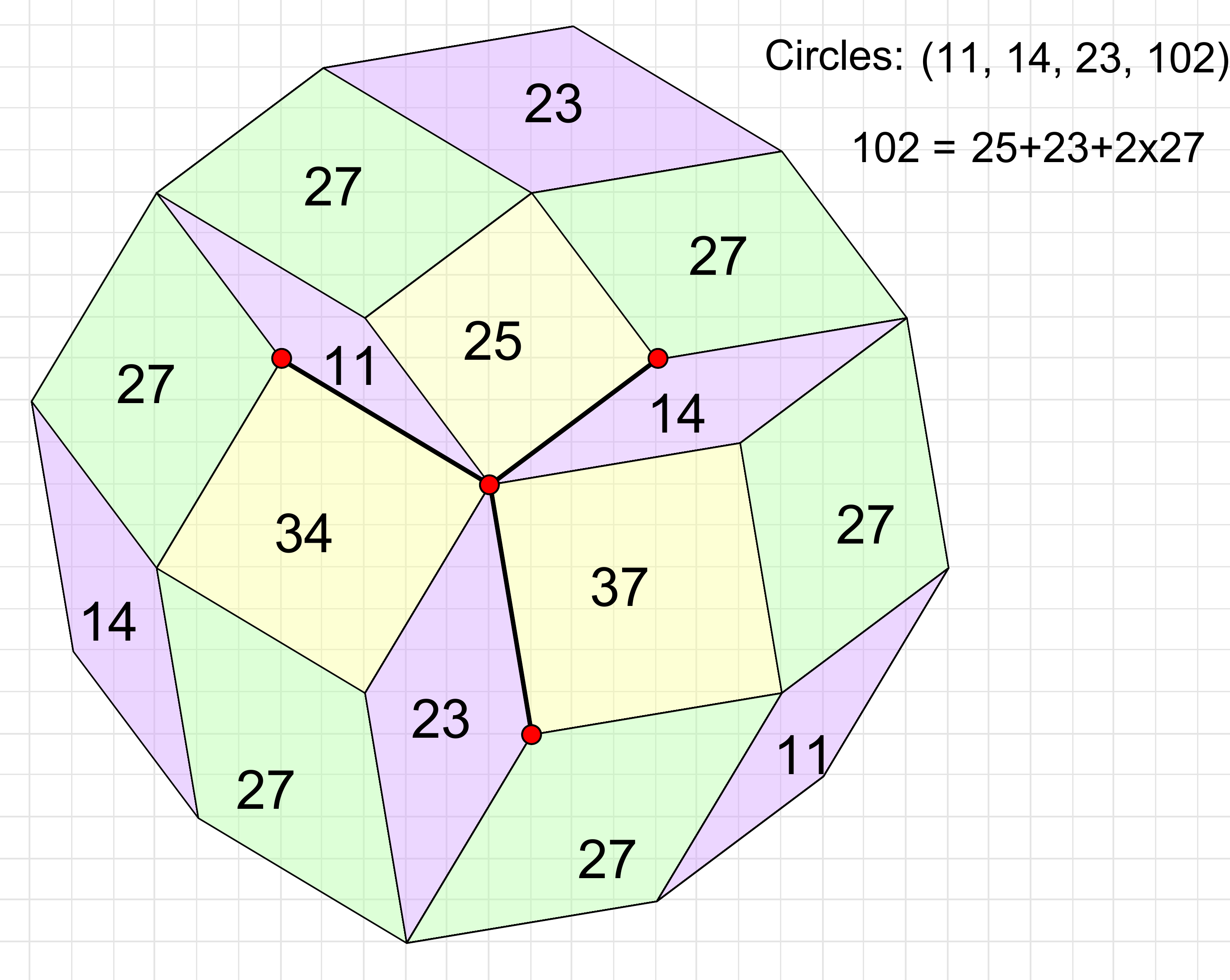}
\caption{ Additional examples: $(3,14,6, 47)$ and $(11,14,23, 102)$}
\label{fig:more}
\end{figure}

\newpage
%---------------------------------------------------------------------------
\section{Spinor space for tangent disks}

In his section we review briefly the notion of spinor description of disks in configurations.
For details,  proofs, and motivation  see \cite{jk-s}.

\subsection{Spinor space}
First,  general notions.  In the following, by the spinor space we will mean 
a two dimensional real vector space, 
which may be also emulated with 1-dimensional complex space (Argand plane), 
$\mathbb T = \mathbb R^2 \cong \mathbb C$.
Typical vectors:
$$
\mathbf a = \begin{bmatrix}  x\\ y  \end{bmatrix} % \ = \ x+yi
\qquad
\mathbf b = \begin{bmatrix}  x' \\ y'  \end{bmatrix}  %\ = \ x'+y' i
\qquad
$$
The space is equipped with two structures,
the ``dot product'' and the ``cross-product'', both with values in real numbers:
$$
\begin{array}{rclll} 
\hbox{inner product:}           &  \mathbf a,\mathbf b  \quad \mapsto\quad   &\mathbf a\bigcdot \mathbf b =  xx' + yy'    \\[5pt]
\hbox{symplectic product:}   &  \mathbf a,\mathbf b  \quad \mapsto\quad   &\mathbf a\times\mathbf  b =  xy' - x'y      \\[5pt]
%\hbox{conjugation:}             &  a \quad \mapsto \quad ai \\[5pt]
%\hbox{relation:}             &  a\times b =   \\[5pt]
\end{array}
$$
We also define `'symplectic conjugation''
$$
\mathbf a^\star \ = \  
         \begin{bmatrix}  x\\ y  \end{bmatrix}^\star \ = \  
         \begin{bmatrix}  -y\\ \phantom{-}x  \end{bmatrix} \ = \ 
         \begin{bmatrix}  0&\!\!\!\!-\!1\\ 1&0  \end{bmatrix}\begin{bmatrix}  x\\ y  \end{bmatrix}
$$
The two structures are related
$$
         \mathbf a\times \mathbf b \ = \  \mathbf a^\star\bigcdot \mathbf b  
$$
There are other identities that are readily implied, e.g., 
$$
 (\mathbf a^\star)^\star = - \mathbf a\,, 
\quad \mathbf a\bigcdot \mathbf b = \mathbf a\times \mathbf b^\star\,, 
\quad \mathbf a^\star \bigcdot \mathbf b^\star = \mathbf a\bigcdot \mathbf b \,,
\quad \mathbf a^\star \times \mathbf b^\star =\mathbf  a\times\mathbf  b \,,
\quad \mathbf a^\star \times \mathbf b = \mathbf b^\star\times \mathbf a \,.
$$
The squares are:
$$
 \mathbf  a\bigcdot \mathbf a =\|\mathbf a\|^2  =x^2 + y^2\,,
\qquad  \mathbf a\times \mathbf a = 0\,,
\qquad \mathbf a\bigcdot \mathbf a^\star = 0\,.
$$
{\bf Interpretation via complex numbers.}
When the spinor space $\mathbb T$ is represented by complex numbers,
the above formulas obtain the following forms.
If
$$
\mathbf a = \begin{bmatrix}  x\\ y  \end{bmatrix}  \ \leftrightarrow \ a= x+yi
\qquad\hbox{and}\qquad
\mathbf b = \begin{bmatrix}  x' \\ y'  \end{bmatrix}  \ \leftrightarrow \ b=  x'+y' i
$$
then the structures are expressed as follows:
$$
\begin{array}{rrlll} 
\hbox{inner product:}           &  a,b  \quad \mapsto\quad   &a\bigcdot b = \frac{1}{2} (\bar a b   + a \bar b)     \\[5pt]
\hbox{symplectic product:}   &  a,b  \quad \mapsto\quad   &a\times b = \frac{1}{2i} (\bar a b   - a \bar b)     \\[5pt]
\hbox{conjugation:}             &  a \quad \mapsto \quad  &a^\star = ai \\[5pt]
%\hbox{relation:}             &  a\times b =   \\[5pt]
\end{array}
$$
Note that by ``conjugation'' we mean ``symplectic conjugation'' (denoted by star).
Not to be confused with ``complex conjugation'',
always called by its full name.
The two descriptions, vector and complex, will be used interchangeably.

%-----------------------------------------------------------------------
\subsection{Spinors and Descartes}

In this section we review the basic facts concerning the ``spinor structure'' of Descartes configurations.
The idea of tangency spinor was introduced in \cite{jk-c}.  For details see \cite{jk-s}.

~\\
{\bf Definition:}  Let $\mathbf A$ and $\mathbf B$ be an ordered pair of mutually tangent disks 
of radii  $r_A$ and $r_B$ and centered at $C_A$ and $C_B$, respectively, 
in a plane identified with complex numbers, $\mathbb C\cong \mathbb R^2$.
Interpret the vector joining the centers as a complex number $z= (C_AC_B)$.
The {\bf tangency spinor} of the two disks is a square root defined 
$$
\mathbf u \  = \ \pm \sqrt{\frac{z}{r_Ar_B}} \ \  \in \mathbb C \cong \mathbb T
\qquad \hbox{and} \qquad 
\mathbf u=\begin{bmatrix}\hbox{Re}\, u\\ \hbox{Im}\, u \end{bmatrix}
$$
We shall view it as a vector of the spinor space $\mathbb T$ discussed in the previous subsection.
The spinor is defined up to a sign since $(-u)^2 = u^2$.
It  depends on the order of disks:  
 if $u$ is a spinor for $(AB)$,  then the spinor for  $(BA)$  is $u^\star = iu$  (again, up to sign).
\\

\begin{figure}[H]
\centering
\includegraphics[scale=.85]{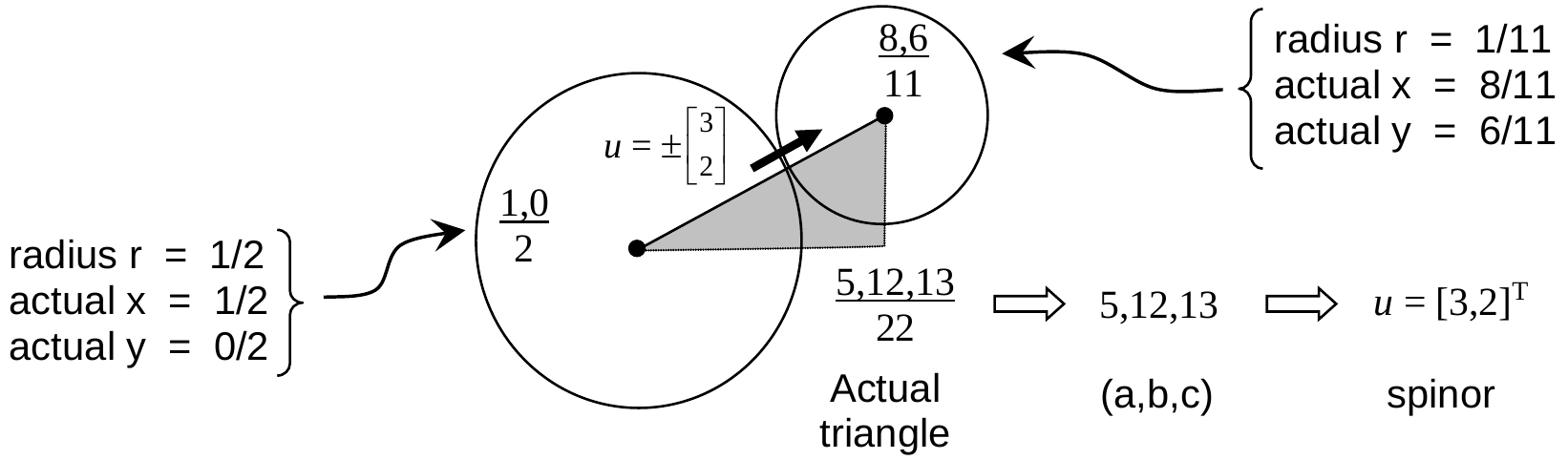}
\caption{\small From two tangent circles to a spinor (not to scale)}
\label{fig:process}
\end{figure}

The geometric interpretation and motivation follows.
Every disk in the Cartesian plane may be given a {\bf symbol},
a fraction-like label that encodes the size and position of the disk:  
the curvature  is indicated in the denominator 
while the positions of the centers may be read off by interpreting the symbol as a pair of fractions
\cite{jk-m}.
%$$
%\hbox{\sf symbol:\  } 
%       \frac{3,\;4}{6}  
%                         \qquad \Longrightarrow \qquad 
%                                     \hbox{\sf radius\  }  r = \frac{1}{6}
%                                      ,\quad  
%                                      \hbox{\sf center } = (x,\,y)  =   \left(\frac{3}{6},\, \frac{4}{6}\right)
%                                                                    =   \left(\frac{1}{2},\, \frac{2}{3}\right)
%$$
$$
\hbox{\sf symbol:\  } 
       \frac{\dot x,\;\dot y}{\beta}  
                         \qquad \Longrightarrow \qquad 
       \begin{cases}     \hbox{\sf radius:\  } &  r = \frac{1}{\beta} \\
                                  \hbox{\sf center:\ }  & (x,\,y)  =   \left(\frac{\dot x}{\beta},\, \frac{\dot y}{\beta}\right) 
       \end{cases}
$$
The numerator, called the {\bf reduced coordinates} of the a disk's center, 
is denoted by dotted letters $(\dot x, \, \dot y)= (x/r,\, y/r)$.
Unbounded disks extending outside a circle are given negative radius and curvature.
Two tangent disks in a plane  define a triangle
with sides as follows:
%$$
\begin{equation}
\label{eq:product}
     \frac{\dot x_1, \;  \dot y_1}{ \beta_1}
\ \Join  \
     \frac{\dot x_2,\; \dot y_2}{ \beta_2}
\qquad \mapsto \qquad
       \left[\begin{matrix}
                     a \\
                     b \\
                    c
     \end{matrix}\right]
\equiv 
  \left[\begin{matrix}
                     \beta_1 \dot x_2 - \beta_2 \dot x_1 \\
                     \beta_1 \dot y_2 - \beta_2 \dot y_1 \\
                    \beta_1+\beta_2
     \end{matrix}\right]
\end{equation}
%$$
where $a^2+b^2=c^2$ 
(see Figure \ref{fig:process}).
The actual size of the triangle in the plane is scaled down by the factor of $\beta_1\beta_2$
(gray triangles in Figure~\ref{fig:process}).
The symbols in some disk packing are integral, then so are the triples $(a,b,c)$.  
Recall that Pythagorean triangles admit Euclidean parameters 
that determine them via the following prescription:
$$
  \mathbf u=\vvec mn \quad \to \quad (a,b,c) = (m^2-n^2,\ 2mn, \ m^2+n^2)
$$
(see, e.g., \cite{Sie, T-T}).
As explained in  \cite{jk-s}, Euclidean parameters can be viewed as a {\it spinor},
a vector of $\mathbb T \cong \mathbf u\in \mathbb R^2$.
Equivalently, viewing the spinor as a complex number $u\in\mathbb C\cong\mathbb R^2$
%via identification $[m, n]^T \equiv m+ni$. 
the above relations is defined by squaring:
$$
            u=
            %\vvec{m}{n} =
            m+ni \quad \to \quad u^2  = a+bi   = (m^2-n^2) + 2mn\, i 
$$
with $c=|u^2|=m^2+n^2$.  
We extend  this map to arbitrary oriented triangles, not necessarily integer.
\\

\newpage

The emergence of the {\bf tangency spinor} for a pair of tangent disks  
is summarized in Figure~\ref{fig:process}.
In graphical representation we shall mark a spinor by an arrow that indicates 
the order of circles, and will label it by its matrix value.
%
%Recall that the spinor is defined up to a sign, since $(-u)^2=u^2$.

%================= FIGURE 2.1

Below, we state the main properties; for proofs see \cite{jk-s}.
The capital letters will denote both circles and their curvatures.
\\
\\

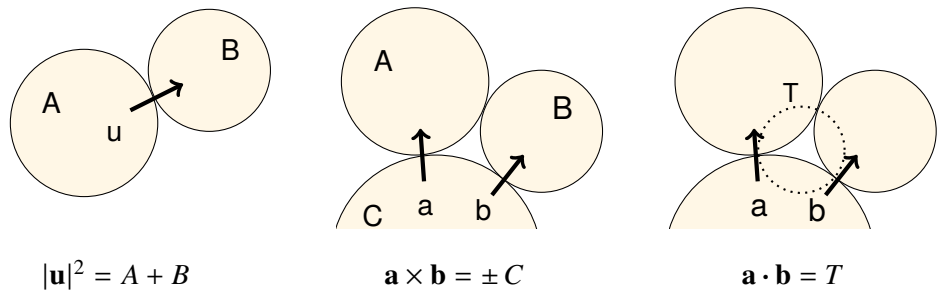
\begin{figure}[H]
%============================ ==================
\centering

%------------ uno
\begin{tikzpicture}[scale=1.4]
\draw[fill=gold!10]  (-.19, 0) circle (.7);
\draw[fill=gold!10]  (1, 0.5) circle (.58);

\draw [->, ultra thick, color=black] (.25, .12) -- (.75, .37);

\node  [scale=1.1]  at (-.5, .2) {\sf  A}; 
\node  [scale=1.1]  at (1.2, .7) {\sf  B}; 

\node  [scale=1.1]  at (.1, -.1) {\sf  u}; 
\draw [color=white, fill=white, opacity=\nofigures, draw opacity=\nofigures] (-.9,-1) rectangle (1.5,1);    
\end{tikzpicture}
\qquad
\begin{tikzpicture}[scale=1.4]
\clip (-.95, .3) rectangle  (1.6,2.45);

\draw[fill=gold!10]  (0,0) circle (1);
\draw[fill=gold!10]  (-.2, 1.7) circle (.7);
\draw[fill=gold!10]  (1, 1.226) circle (.58);

\draw [->, ultra thick, color=black] (-.115, .75) -- (-.15,1.25);
\draw [->, ultra thick, color=black] (.53, .62) -- (.83,1.0);

\node  [scale=1.1]  at (-.6, .43) {\sf  C}; %   {$\frac{\a,\b}{\c,\d}$}
\node  [scale=1.1]  at (-.5, 1.9) {\sf  A}; %   {$\frac{\a,\b}{\c,\d}$}
\node  [scale=1.2]  at (1.2, 1.45) {\sf  B}; %   {$\frac{\a,\b}{\c,\d}$}
\

\node  [scale=1.1]  at (-.1, .5) {\sf  a}; 
\node  [scale=1.1]  at (.45, .48) {\sf  b}; 
\draw [color=white, fill=white, opacity=\nofigures, draw opacity=\nofigures] (-.9,.3) rectangle (1.6,2.4);    
\end{tikzpicture}
\qquad
\begin{tikzpicture}[scale=1.4]
\clip (-.95, .3) rectangle  (1.6,2.45);

\draw[fill=gold!10]  (0,0) circle (1);
\draw[fill=gold!10]  (-.2, 1.7) circle (.7);
\draw[fill=gold!10]  (1, 1.226) circle (.58);
\draw[thick, dotted]  (.3, 1.05) circle (.41);

\draw [->, ultra thick, color=black] (-.115, .75) -- (-.15,1.25);
\draw [->, ultra thick, color=black] (.53, .62) -- (.83,1.0);

\node  [scale=1.2]  at (-.1, .5) {\sf  a}; 
\node  [scale=1.2]  at (.45, .48) {\sf  b}; 
\node  [scale=1]  at (.2, 1.6) {\sf T}; 
\end{tikzpicture}

\begin{tikzpicture}[scale=1.5]
\node  [scale=1]  at (-3,0) {$|\mathbf u|^2 = A+B $}; 
\node  [scale=1]  at (0, 0) {$\mathbf a\times \mathbf b =\pm \, C$}; 
\node  [scale=1]  at (3, 0) {$\mathbf a\bigcdot \mathbf b = T$};
\node  [scale=1]  at (4, 0) {~};
\node  [scale=1]  at (0,.3) {~};
%\draw [color=white, fill=white, opacity=\nofigures, draw opacity=\nofigures] (-.2,.5) rectangle (1.4,2);    
\end{tikzpicture}

\caption{Spinors and disks.}
\label{fig:3A}
\end{figure}

\noindent
{\bf Proposition 1.} % ---------------------prop   7.3
\label{thm:curvs}
 If $u$ is the tangency spinor for two tangent disks of curvatures $A$ and $B$, respectively,   
(Figure \ref{fig:3A}, left)  then 
\begin{equation}
\label{eq:thm3}
                     |\mathbf u|^2 = A + B
\end{equation}
%
%\end{theorem}

\noindent
{\bf Theorem 2}  [{\bf curvatures from spinors}]    %-------------------------Theorem 7.4  ]:   
\label{thm:curv}
In the system of three mutually tangent circles, the symplectic product of two spinors directed outward from 
(respectively inward into) one of the circles equals (up to sign) its curvature, 
e.g., following notation of Figure \ref{fig:3A}, center:
\begin{equation}
\label{eq:thm4}
          \pm \, C  \ = \ \mathbf a\times \mathbf b \   =   \ \det [\mathbf a|\mathbf b]
\end{equation}

%%%%%%%%%%%%%%%%%%%%%%%%%%%%%%%%%%%%%%%%%%%

\newpage

\noindent
{\bf Theorem 3}  [{\bf mid-circles from spinors}].    %-------------------------Theorem 7.4  ]:   
\label{thm:curv}
In the system of three mutually tangent circles, the dot product of two spinors directed outward
(respectively, inward)  from 
one of the circles equals to the curvature
of the mid-circle $(ABC)$. Following the notation of Figure \ref{fig:3A}, right:
\begin{equation}
\label{eq:thm42}
           T  =  \mathbf a\bigcdot \mathbf b \qquad (\hbox{respectively, } T = - \mathbf a\bigcdot \mathbf b)
\end{equation}
In particular, $\mathbf a\bigcdot \mathbf b = -\mathbf a^\star \bigcdot \mathbf b^\star$.

\begin{figure}[H]
%============================ ==================
\centering

\begin{tikzpicture}[scale=1.4]
\clip (-.95, .3) rectangle  (1.6,2.45);

\draw[fill=gold!10]  (0,0) circle (1);
\draw[fill=gold!10]  (-.2, 1.7) circle (.7);
\draw[fill=gold!10]  (1, 1.226) circle (.58);

\draw [->, ultra thick, color=black] (-.115, .75) -- (-.15,1.25);
\draw [<-, ultra thick, color=black] (.53, .62) -- (.83,1.0);
\draw [->, ultra thick, color=black] (.25, 1.54) -- (.75,1.33);

\node  [scale=1.1]  at (-.1, .5) {\sf  a}; 
\node  [scale=1.1]  at (1, .98) {\sf  c}; 
\node  [scale=1.1]  at (.1, 1.7) {\sf  b}; 
\draw [color=white, fill=white, opacity=\nofigures, draw opacity=\nofigures] (-.9,.3) rectangle (1.6,2.4);    
\end{tikzpicture}
\qquad
\begin{tikzpicture}[scale=1.4]
\clip (-.95, .3) rectangle  (1.6,2.45);
\draw[fill=gold!10]  (0,0) circle (1);
\draw[fill=gold!10]  (-.2, 1.7) circle (.7);
\draw[fill=gold!10]  (1, 1.226) circle (.58);
\draw [fill=gold!10] (.32, 1.07) circle (.1167);

\draw [->, ultra thick, color=black] (.115, .6) -- (.24,.9);
\node  [scale=1.1]  at (-.1, .5) {\sf  a}; 

\draw [->, ultra thick, color=black] (.83, 1.22) -- (.53,1.12);
\node  [scale=1.1]  at (1, 1.2) {\sf  c}; %   {$\frac{\a,\b}{\c,\d}$}

\draw [->, ultra thick, color=black] (.03, 1.45) -- (.23,1.2);
\node  [scale=1.1]  at (0, 1.7) {\sf  b}; 
\draw [color=white, fill=white, opacity=\nofigures, draw opacity=\nofigures] (-.9,.3) rectangle (1.6,2.4);    
\end{tikzpicture}
\qquad
\begin{tikzpicture}[scale=1.4]
\clip (-3, -.08) rectangle  (0.1, 2.2);

\draw[fill=gold!10]  (0,-5) circle (6);
\draw[fill=gold!10]  (-2, 1.4) circle (.7);
\draw[fill=gold!10]  (-.8, 1.455) circle (.5);
\draw[fill=gold!10]  (-1.257, 1.01) circle (.136);

\draw [->, thick, color=black] (-1.6, .22) -- (-1.75, .67);
\node  [scale=1.1]  at (-1.55, .01) {\sf  a}; 

\draw [->, thick, color=black] (-.65, .37) -- (-.73, .89);
\node  [scale=1.1]  at (-.58, .14) {\sf  b}; 

\draw [->, thick, color=black] (-1.12, .51) -- (-1.2, .82);
\node  [scale=1]  at (-1.1, .28) {\sf  a+b}; 
\end{tikzpicture}

\begin{tikzpicture}[scale=1.5]
\node  [scale=1]  at (-3,0) {$\mathbf a+\mathbf b+\mathbf c \ =\ 0$}; 
\node  [scale=1]  at (0, 0) {$\mathbf a+\mathbf b+\mathbf c \ =\ 0$}; 
\node  [scale=1]  at (3, 0) {$\mathbf a+\mathbf b \ =\ \mathbf c$};
\node  [scale=1]  at (4, 0) {~};
\node  [scale=1]  at (0,.3) {~};
%\draw [color=white, fill=white, opacity=\nofigures, draw opacity=\nofigures] (-3,.-.08) rectangle (.1, 2.2);    
\end{tikzpicture}

\caption{Spinor fundamental theorem}
\label{fig:3B}
\end{figure}
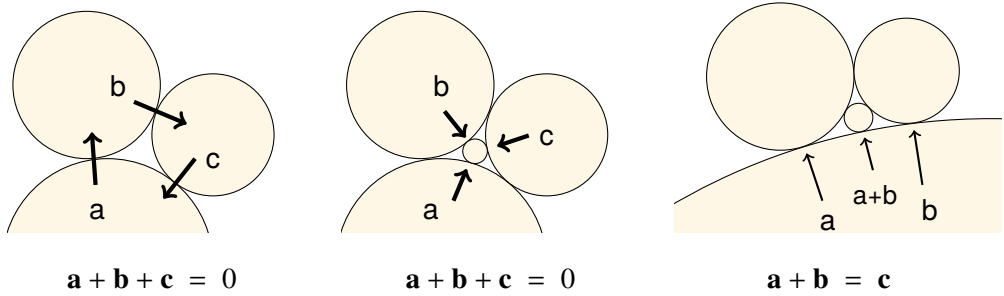
%%%%%%%%%%%%%%%%%%%%%%%%%%%%%%%%%%%%%%%%%%%

\noindent
{\bf Theorem 4}  [{\bf spinor curl}]. %----------------------------
\label{thm:curl}
The signs of the three tangency spinors between be three mutually tangent circles
(Figure \ref{fig:3B}, left)  
may be chosen so that 
\begin{equation}
\label{eq:curl}
\mathbf a  +  \mathbf b  +  \mathbf c  =  \mathbf 0  \qquad\quad        
 [ \hbox{\rm ``curl}\ \mathbf u = 0 \hbox{''}]   %\bm u
\end{equation}
%
%\end{theorem}

The  label in the above equation is to be understood figuratively meaning, and ``$\mathbf  u$'' in it stands for the triple of 
spinor $(\mathbf a,\mathbf b, \mathbf c)$.

~

\noindent
%-------------- Theorem 7.4 : 
{\bf Theorem 5.}
%\label{thm:div}
Let $A$, $B$, $C$, and $D$ be four circles in a Descartes configuration.  
\\[5pt]
%%%%%%%%%%%   - -----   -\begin{addmargin}[1em]{0em}
{\bf [A. Vanishing divergence]:}  
If $a$, $b$ and $c$ are tangency spinors for pairs $AD$, $BD$ and $CD$  
(see Figure \ref{fig:3B}, center),
then their signs may be chosen so that 
\begin{equation}
\label{eq:thm5a}
 					\mathbf a + \mathbf b + \mathbf c = 0 	\qquad	[\hbox{``\rm div}\, \mathbf u = 0\hbox{''}]                                   
\end{equation}
The same property holds for the outward oriented spinors.
\\[5pt]
{\bf [B. Additivity]:}  If $a$ and $b$ are spinors of tangency for pairs $CA$ and $CB$ 
(see  Figure \ref{fig:3B}, right), 
then there is a choice of signs so that the sum 
\begin{equation}
\label{eq:thm5b}
                                      \mathbf c = \mathbf a +\mathbf  b
\end{equation}
is the tangency  spinor of $CD$.  

%----------------------------\end{addmargin}
%\end{theorem}

%%================= FIGURE 2.1
%\begin{figure}[h]
%\centering
%\includegraphics[scale=.7]{S75}
%\caption{\small (a) vanishing divergence, (b) spinor addition}
%\label{fig:thmdiv}
%\end{figure}
%

\newpage

%----------------------------------------------------------------------
\section{Proof of the correspondence}

Here is the theorem that summarizes the observations of Section 1.
%The general statements are accompanied with numerical values as in Section~1. 

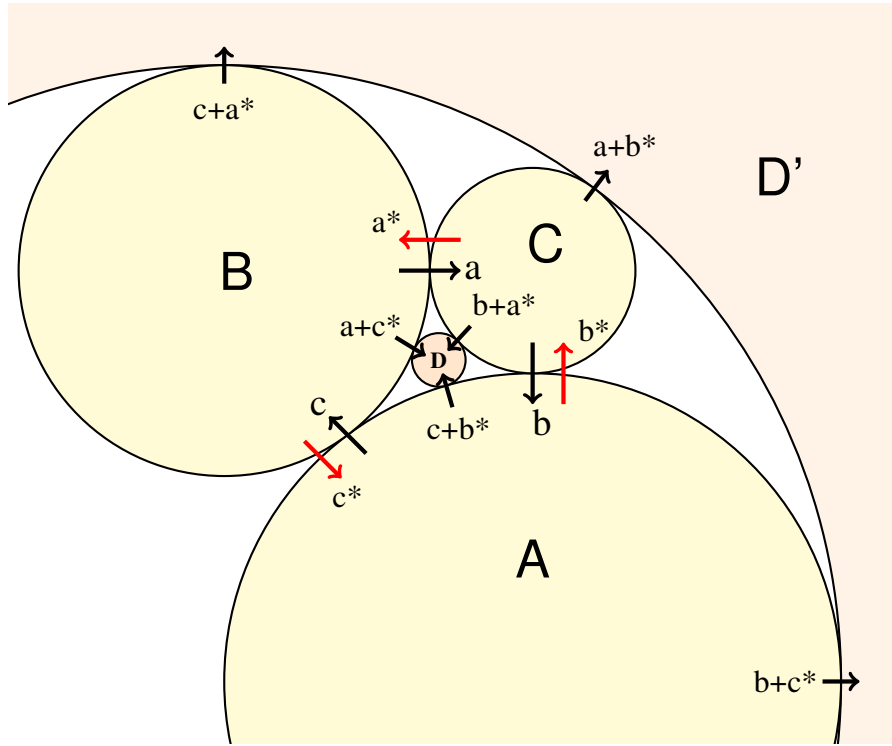
\begin{figure}[H]
\centering
\begin{tikzpicture}[scale=8.2, rotate=0]
\clip (-.35,-.1) rectangle (1.1,  1.1);
%\draw (0,0) circle (1);
%\draw [thick] (0,0) circle (1);
%\draw [fill=blue, opacity=.125, thick] (0,0) circle (1);
% A  1, 0, 2                  sums  12,  18,  34
% B  0, 2, 3
% C 3, 4, 6
%  D  8, 12 23 

\draw [yellow!20, fill=orange!10]   (-.55, -.5) rectangle (1.5, 1.5);

\draw [thick, fill=white] (0,0) circle (1)
node  [scale=2]  at (.9, .82) {\sf  D'}; %   {$\frac{\a,\b}{\c,\d}$}

\draw [thick, fill=yellow!20] (1/2, 0) circle (1/2)
node  [scale=2]  at (.5, .2) {\sf  A}; %   {$\frac{\a,\b}{\c,\d}$}

\draw [thick, fill=yellow!20] (0,2/3) circle (1/3)
node  [scale=2]  at (1/50,  2/3) {\sf B}; %   {$\frac{\a,\b}{\c,\d}$}

\draw [thick, fill=yellow!20] (3/6, 4/6) circle (1/6)
node  [scale=2]  at (3/6 +.021,  4/6 +.045) {\sf C}; %   {$\frac{\a,\b}{\c,\d}$}

% (A+B+C-D')/2   
\draw [thick, fill=orange!20] (8/23, 12/23) circle (1/23);
\node  [scale=.79]  at (8/23,12/23) {\sf\bf D}; %   {$\frac{\a,\b}{\c,\d}$}

%% (A+B+C-D)/2   2, 3, 6
%\draw [thick, dotted] (2/6, 3/6) circle (1/6);
%% (A+B+D-C)/2     3, 5, 11
%\draw [thick, dotted, fill=yellow!50, fill opacity=.70] (3/11, 5/11) circle (1/11);
%% (A+C+D-B)/2     6, 7, 14
%\draw [thick, dotted, fill=yellow!50, fill opacity=.70] (6/14, 7/14) circle (1/14);
%% (B+C+D-A)/2     5, 9, 15
%\draw [thick, dotted, fill=yellow!50, fill opacity=.70] (5/15, 9/15) circle (1/15);
%
%\draw node  [scale=1.1]  at (.4,.37) {\color{red} \sf  6}; %   {$\frac{\a,\b}{\c,\d}$}
%\draw node  [scale=1.1]  at (.24,.49) {\color{red} \sf 11}; %   {$\frac{\a,\b}{\c,\d}$}
%\draw node  [scale=1.1]  at (.44,.54) {\color{red} \sf  14}; %   {$\frac{\a,\b}{\c,\d}$}
%\draw node  [scale=1.1]  at (.30,.61) {\color{red} \sf  15}; %   {$\frac{\a,\b}{\c,\d}$}

\draw [->, ultra thick, color=black] (1/3-.05, 2/3) -- (1/3 +.05, 2/3);
\draw node at (1/3+.07, 2/3)  {\Large a};
\draw [->, ultra thick, color=black] (1/2, 1/2 + .05) -- (1/2, 1/2 -.05);
\draw node at (1/2+.015, 1/2 - .079) {\Large b};
\draw [->, ultra thick, color=black] (1/5+.03, 2/5 - .03) -- (1/5-.03, 2/5 +.03);
\draw node at (1/5-.047, 2/5 + .03+ .017) {\Large c};

\draw [->, ultra thick, color=red] (1/3+.05, 2/3+.05) -- (1/3 -.05, 2/3+.05);
\draw node at (1/3-.07, 2/3+.08)  {\large a*};
\draw [->, ultra thick, color=red] (1/2 + .05, 1/2 - .05) -- (1/2+.05, 1/2 +.05);
\draw node at (1/2+.10, 1/2+.07)  {\large b*};
\draw [->, ultra thick, color=red] (1/5-.03-.04, 2/5 + .03-.04) -- (1/5+.03-.04, 2/5 -.03-.04);
\draw node at (1/5, 2/5-.1)  {\large c*};

\draw [->, ultra thick, color=black] (9/25+.01, 12/25-.035) -- (9/25-.006 , 12/25+.02);
\draw node at (9/25+.02, 12/25-.07)  {\large c+b*};
\draw [->, ultra thick, color=black] (8/26-.03, 14/26+.02) -- (8/26+.019 , 14/26-.01);
\draw node at (8/26-.07, 14/26+.04)  {\large a+c*};
\draw [->, ultra thick, color=black] (11/29+.02, 16/29+.025) -- (11/29-.017 , 16/29-.015);
\draw node at (11/29 + .075 , 16/29 +.06)  {\large b+a*};

\draw [->, ultra thick, color=black] (0, 1-.03) -- (0 , 1+.03);
\draw node at (0, 1-.07)  {\large c+a*};
\draw [->, ultra thick, color=black] (1-.03, 0) -- (1+.03 , 0);
\draw node at (1-.09, 0)  {\large b+c*};
\draw [->, ultra thick, color=black] (3/5-.015, 4/5-.02) -- (3/5+.0225, 4/5+.03);
\draw node at (3/5+.05, 4/5+.07)  {\large a+b*};

\draw [color=white, fill=white, opacity=\nofigures, draw opacity=\nofigures] (-1,-1) rectangle (1,1);     %<<<<<<<<<<<<<<<<<<<<<<<<<<<
\end{tikzpicture}
\caption{Distribution of spinors in Descartes configuration (A,B,C,D) and (A,B,C,D')}
\label{fig:distribution} 
\end{figure}

\def\Curv{\hbox{\sf Curv}\/}
\def\Area{\hbox{\sf Area}\/}

%----------------------
\noindent
{\bf Theorem 6:}  Suppose $A$, $B$, $C$, and $D$ are disks in Descartes configuration. 
Let $a$, $b$, $c$ are tangency spinors joining $A$, $B$, and $C$, as shown in Figure~\ref{fig:distribution}.  
Then the ``three-spinor tessellation''  built on these spinors 
determines curvatures of disks and circles as follows:
\begin{enumerate}
\item
The curvatures $A$, $B$ and $C$ are equal to the areas of the three red tiles.
The fourth inscribed disk (solution to Descartes problem) has curvature equal to the butterfly 
(see Fig.~\ref{fig:3butterflies}):  
$$
\begin{array}{rl}
  \Curv(D)  &=  \Area( \, 3\;  \hbox{red} + \hbox{two green}\,)\\
                &=  \Area( \, \hbox{yellow} + \hbox{aligned red} + \hbox{two green} \, )  \,.
\end{array}
$$
Equivalently,one may consider the ring of the three red (one of them light red) and two green tiles. 

\item
All green tiles have the same area, equal to the curvature of the mid-circle $(ABC)$
that passes through the points of tangency of disks A,B,C.
This, and the other mid-circles have curvatures:
$$
\begin{array}{rl}
\Curv (ABC) &= \Area( \, \hbox{any green}  \, )  \\
\Curv (ABD) &= \Area( \,  \hbox{yellow} (A+B)  + \hbox{green} \,  )\\
\Curv (BCD) &= \Area( \,  \hbox{yellow} (B+C) + \hbox{green} \,  )\\
\Curv (CAD) &= \Area( \,  \hbox{yellow} (C+A) + \hbox{green} \,  )
\end{array}
$$
\item
The second solution to the Descartes problem for disks $A$, $B$, and $C$, is the disk of curvature 
$$
\begin{array}{rl}
  \Curv (D') & = \Area( \,  3 \; \hbox{red} - \hbox{two green} \,) \\
                 & = \Area( \, \hbox{yellow} + \hbox{aligned red}  - \hbox{two green} \, ) \,.
\end{array}
$$
%
%here 9+2-2 6 = 8+3-2 6 = 5+6 - 2.6 = -1

The mid-circles for this Descartes configuration are 
$$
\begin{array}{rl}
 \Curv (ABD')  &= \Area( \, \hbox{yellow}(A+B) - \hbox{green} \, )\\
 \Curv (BCD')  &= \Area( \, \hbox{yellow}(B+C) - \hbox{green} \, )\\
  \Curv (CAD') &= \Area( \, \hbox{yellow}(C+A) - \hbox{green} )
  \end{array}
  $$
\end{enumerate}

%===============================================================
\begin{figure}[t]
\centering
\begin{tikzpicture}[scale=.75, rotate=0, shift={(0,0cm)}]  %was 7
\draw (-5,-5) grid (5,5);
%%\clip (-1.25,-2.1) rectangle (1.1,2.2);
%\draw (-\ile,0) -- (\ile,0);

\def\Ax{3}; \def\Ay{0};
\def\Bx{-1}; \def\By{2};
\def\Cx{-2}; \def\Cy{-2};
\def\ax{0-\Ay}; \def\ay{\Ax}; \def\bx{0-\By}; \def\by{\Bx}; \def\cx{0-\Cy}; \def\cy{\Cx};

%-----------------yellow
\draw [thick, fill=yellow!50, fill opacity=.79]  (0,0) -- (\Ax,\Ay) --  (\Ax+\ax, \Ay + \ay)  -- (\ax,\ay) -- cycle;
\draw [thick, fill=yellow!50, fill opacity=.79]  (0,0) -- (\Bx, \By) --  (\Bx+\bx, \By + \by)   -- (\bx,\by) -- cycle;
\draw [thick, fill=yellow!50, fill opacity=.79]  (0,0) -- (\Cx,\Cy) --  (\Cx+\cx,\Cy+\cy) -- (\cx,\cy) -- cycle;
%------------------red
\draw [thick, fill=red!80, fill opacity=.79]  (0,0) -- (\ax,\ay) -- (\ax+\Bx,\ay+\By) -- (\Bx,\By) -- cycle;
\draw [thick, fill=red!80, fill opacity=.79]  (0,0) -- (\bx,\by) -- (\bx+\Cx,\by+\Cy) -- (\Cx,\Cy) -- cycle;
\draw [thick, fill=red!80, fill opacity=.79]  (0,0) -- (\cx,\cy) -- (\cx+\Ax,\cy+\Ay) -- (\Ax,\Ay) -- cycle;
%-----------------green
%aroun yellow 1, left, right
\draw [thick, fill=green!50, fill opacity=.79]  (\ax,\ay) -- (\Ax+\ax,\Ay+\ay) 
                                                                            -- (\Ax+\ax+\Bx,\Ay+\ay+\By) -- (\ax+\Bx,\ay+\By) -- cycle;
\draw [thick, fill=green!50, fill opacity=.79]  (\Ax,\Ay) -- (\Ax+\ax,\Ay+\ay) 
                                                                            -- (\Ax+\ax+\cx,\Ay+\ay+\cy) -- (\Ax+\cx,\Ay+\cy) -- cycle;                                                                   
%aroun yellow 2, left, right
\draw [thick, fill=green!50, fill opacity=.79]  (\bx,\by) -- (\Bx+\bx,\By+\by) 
                                                                            -- (\Bx+\bx+\Cx,\By+\by+\Cy) -- (\bx+\Cx,\by+\Cy) -- cycle;                                                                            
\draw [thick, fill=green!50, fill opacity=.79]  (\Bx,\By) -- (\Bx+\bx,\By+\by) 
                                                                            -- (\Bx+\bx+\ax,\By+\by+\ay) -- (\Bx+\ax,\By+\ay) -- cycle;                                                                   
%aroun yellow 3, left, right
\draw [thick, fill=green!50, fill opacity=.79]  (\cx,\cy) -- (\Cx+\cx,\Cy+\cy) 
                                                                            -- (\Cx+\cx+\Ax,\Cy+\cy+\Ay) -- (\cx+\Ax,\cy+\Ay) -- cycle;                                                                            
\draw [thick, fill=green!50, fill opacity=.79]  (\Cx,\Cy) -- (\Cx+\cx,\Cy+\cy) 
                                                                            -- (\Cx+\cx+\bx,\Cy+\cy+\by) -- (\Cx+\bx,\Cy+\by) -- cycle;                                                                                
                                                                          
%--------------------lightred
\draw [thick, fill=red!50, fill opacity=.79]   (\Ax+\ax,\Ay+\ay) -- (\Ax+\ax+\cx,\Ay+\ay+\cy)
                                                        -- (\Ax+\ax+\Bx+\cx,\Ay+\ay+\By+\cy) -- (\Ax+\ax+\Bx,\Ay+\ay+\By)   -- cycle;                                                                            
\draw [thick, fill=red!50, fill opacity=.79]   (\Bx+\bx,\By+\by) -- (\Bx+\bx+\ax,\By+\by+\ay)
                                                        -- (\Bx+\bx+\Cx+\ax,\By+\by+\Cy+\ay) -- (\Bx+\bx+\Cx,\By+\by+\Cy)   -- cycle;                                                                            
\draw [thick, fill=red!50, fill opacity=.79]   (\Cx+\cx,\cy+\cy) -- (\Cx+\cx+\bx,\Cy+\cy+\by)
                                                        -- (\Cx+\cx+\Ax+\bx,\Cy+\cy+\Ay+\by) -- (\Cx+\cx+\Ax,\Cy+\cy+\Ay)   -- cycle;                                                                            

%lines and dots
\draw [fill=black] (0,0) circle  (.2);
\draw [fill=black]  (\Ax,\Ay) circle  (.2);
\draw [fill=black]  (\Bx,\By) circle  (.2);
\draw [fill=black]  (\Cx,\Cy) circle  (.2);
\draw [line width=2pt] (0,0) -- (\Ax,\Ay);
\draw [line width=2pt] (0,0) -- (\Bx,\By);
\draw [line width=2pt] (0,0) -- (\Cx,\Cy);

\node at (1.5,.22) [scale=1.2, color=black] {$\mathbf a$};
\node at (-.42, 1.49) [scale=1.2, color=black] {$\mathbf b$};
\node at (-1.1,-1.5) [scale=1.2, color=black] {$\mathbf c$};

\node at (1.5,1.5) [scale=1.5, color=black] {$|\mathbf a|^2$};
\node at (-1.5,.5) [scale=1.5, color=black] {$|\mathbf b|^2$};
\node at (0,-2) [scale=1.5, color=black] {$|\mathbf c|^2$};

\node at (1,4) [scale=1.4, color=black] {$\mathbf a\!\times\! \mathbf b$};
\node at (-3.5,-1) [scale=1.4, color=black] {$\mathbf b\!\times\! \mathbf c$};
\node at (2.5,-3) [scale=1.4, color=black] {$\mathbf b\!\times\! \mathbf c$};

\node at (4.1,.55) [scale=1.4, color=black] {$\mathbf c^\star\!\!\times \!  \mathbf a^\star$};
\node at (-2+.02,3) [scale=1.4, color=black] {$\mathbf a^\star\!\!\times\! \mathbf b ^\star$};
\node at (-2,-3.5) [scale=1.4, color=black] {$\mathbf b^\star\!\times \mathbf c^\star$};

\node at (2.5,-1) [scale=1.4, color=black] {$\mathbf c^\star\!\!\times\! \mathbf a$};
\node at (-0.5,2.58) [scale=1, color=black] {$\mathbf a^\star\!\!\!\!\times\!\! \mathbf b$};
\node at (-2,-1.45) [scale=1, rotate=37, color=black] {$\mathbf b^\star\!\!\times\!\mathbf c$}; 

\node at (3.5,3.05) [scale=1, rotate=-42, color=black] {$\mathbf c^\star\!\!\times\! \mathbf b$};
\node at (-4,1.5) [scale=1.4,  color=black] {$\mathbf a^\star\!\!\times\!\mathbf  c$};
\node at (0.42,-4.5) [scale=1.4, color=black] {$\mathbf b^\star\!\!\times\! \mathbf a$};

%\draw [fill=gold!10] (.75,.75/2) circle  (.75/2);
\draw [white, fill=white, fill opacity=\nofigures, draw opacity=\nofigures] (-5,-5) rectangle (5.1,5);     %<<<<<<<<<<<<<<<<<<<<<<<<<<<
\end{tikzpicture}
\quad
\begin{tikzpicture}[scale=.75, rotate=0, shift={(0,0cm)}]  %was 7
\draw (-5,-5) grid (5,5);
%%\clip (-1.25,-2.1) rectangle (1.1,2.2);
%\draw (-\ile,0) -- (\ile,0);

\def\Ax{3}; \def\Ay{0};
\def\Bx{-1}; \def\By{2};
\def\Cx{-2}; \def\Cy{-2};
\def\ax{0-\Ay}; \def\ay{\Ax}; \def\bx{0-\By}; \def\by{\Bx}; \def\cx{0-\Cy}; \def\cy{\Cx};

%-----------------yellow
\draw [thick, fill=yellow!50, fill opacity=.79]  (0,0) -- (\Ax,\Ay) --  (\Ax+\ax, \Ay + \ay)  -- (\ax,\ay) -- cycle;
\draw [thick, fill=yellow!50, fill opacity=.79]  (0,0) -- (\Bx, \By) --  (\Bx+\bx, \By + \by)   -- (\bx,\by) -- cycle;
\draw [thick, fill=yellow!50, fill opacity=.79]  (0,0) -- (\Cx,\Cy) --  (\Cx+\cx,\Cy+\cy) -- (\cx,\cy) -- cycle;
%------------------red
\draw [thick, fill=red!80, fill opacity=.79]  (0,0) -- (\ax,\ay) -- (\ax+\Bx,\ay+\By) -- (\Bx,\By) -- cycle;
\draw [thick, fill=red!80, fill opacity=.79]  (0,0) -- (\bx,\by) -- (\bx+\Cx,\by+\Cy) -- (\Cx,\Cy) -- cycle;
\draw [thick, fill=red!80, fill opacity=.79]  (0,0) -- (\cx,\cy) -- (\cx+\Ax,\cy+\Ay) -- (\Ax,\Ay) -- cycle;
%-----------------green
%aroun yellow 1, left, right
\draw [thick, fill=green!50, fill opacity=.79]  (\ax,\ay) -- (\Ax+\ax,\Ay+\ay) 
                                                                            -- (\Ax+\ax+\Bx,\Ay+\ay+\By) -- (\ax+\Bx,\ay+\By) -- cycle;
\draw [thick, fill=green!50, fill opacity=.79]  (\Ax,\Ay) -- (\Ax+\ax,\Ay+\ay) 
                                                                            -- (\Ax+\ax+\cx,\Ay+\ay+\cy) -- (\Ax+\cx,\Ay+\cy) -- cycle;                                                                   
%aroun yellow 2, left, right
\draw [thick, fill=green!50, fill opacity=.79]  (\bx,\by) -- (\Bx+\bx,\By+\by) 
                                                                            -- (\Bx+\bx+\Cx,\By+\by+\Cy) -- (\bx+\Cx,\by+\Cy) -- cycle;                                                                            
\draw [thick, fill=green!50, fill opacity=.79]  (\Bx,\By) -- (\Bx+\bx,\By+\by) 
                                                                            -- (\Bx+\bx+\ax,\By+\by+\ay) -- (\Bx+\ax,\By+\ay) -- cycle;                                                                   
%aroun yellow 3, left, right
\draw [thick, fill=green!50, fill opacity=.79]  (\cx,\cy) -- (\Cx+\cx,\Cy+\cy) 
                                                                            -- (\Cx+\cx+\Ax,\Cy+\cy+\Ay) -- (\cx+\Ax,\cy+\Ay) -- cycle;                                                                            
\draw [thick, fill=green!50, fill opacity=.79]  (\Cx,\Cy) -- (\Cx+\cx,\Cy+\cy) 
                                                                            -- (\Cx+\cx+\bx,\Cy+\cy+\by) -- (\Cx+\bx,\Cy+\by) -- cycle;                                                                                
                                                                          
%--------------------lightred
\draw [thick, fill=red!50, fill opacity=.79]   (\Ax+\ax,\Ay+\ay) -- (\Ax+\ax+\cx,\Ay+\ay+\cy)
                                                        -- (\Ax+\ax+\Bx+\cx,\Ay+\ay+\By+\cy) -- (\Ax+\ax+\Bx,\Ay+\ay+\By)   -- cycle;                                                                            
\draw [thick, fill=red!50, fill opacity=.79]   (\Bx+\bx,\By+\by) -- (\Bx+\bx+\ax,\By+\by+\ay)
                                                        -- (\Bx+\bx+\Cx+\ax,\By+\by+\Cy+\ay) -- (\Bx+\bx+\Cx,\By+\by+\Cy)   -- cycle;                                                                            
\draw [thick, fill=red!50, fill opacity=.79]   (\Cx+\cx,\cy+\cy) -- (\Cx+\cx+\bx,\Cy+\cy+\by)
                                                        -- (\Cx+\cx+\Ax+\bx,\Cy+\cy+\Ay+\by) -- (\Cx+\cx+\Ax,\Cy+\cy+\Ay)   -- cycle;                                                                            

%lines and dots
\draw [fill=black] (0,0) circle  (.2);
\draw [fill=black]  (\Ax,\Ay) circle  (.2);
\draw [fill=black]  (\Bx,\By) circle  (.2);
\draw [fill=black]  (\Cx,\Cy) circle  (.2);
\draw [line width=2pt] (0,0) -- (\Ax,\Ay);
\draw [line width=2pt] (0,0) -- (\Bx,\By);
\draw [line width=2pt] (0,0) -- (\Cx,\Cy);

\node at (1.5,.22) [scale=1.2, color=black] {$\mathbf a$};
\node at (-.5, 1.49) [scale=1.2, color=black] {$\mathbf b$};
\node at (-1.1,-1.5) [scale=1.2, color=black] {$\mathbf c$};

\node at (1.5,1.5) [scale=1.5, color=black] {$B+C$};
\node at (-1.5,.5) [scale=1.5, color=black] {$C+A$};
\node at (0,-2) [scale=1.5, color=black] {$A+B$};

\node at (1,4) [scale=1.4, color=black] {$G$};
\node at (-3.5,-1) [scale=1.4, color=black] {$G$};
\node at (2.5,-3) [scale=1.4, color=black] {$G$};

\node at (4.1,.5) [scale=1.4, color=black] {$G$};
\node at (-2,3) [scale=1.4, color=black] {$G$};
\node at (-2,-3.5) [scale=1.4, color=black] {$G$};

\node at (2.5,-1) [scale=1.4, color=black] {$B$};
\node at (-0.5,2.5) [scale=1.4, color=black] {$C$};
\node at (-2,-1.42) [scale=1.4, color=black] {$A$};

\node at (3.55,3) [scale=1.4, color=black] {$A$};
\node at (-4,1.5) [scale=1.4, color=black] {$B$};
\node at (0.42,-4.5) [scale=1.4, color=black] {$C$};

%\draw [fill=gold!10] (.75,.75/2) circle  (.75/2);
\draw [white, fill=white, fill opacity=\nofigures, draw opacity=\nofigures] (-5,-5) rectangle (5.1,5);     %<<<<<<<<<<<<<<<<<<<<<<<<<<<
\end{tikzpicture}
\caption{Left: Areas of tiles.  Right: Corresponding disk curvatures.}
\label{fig:labeled}
\end{figure}
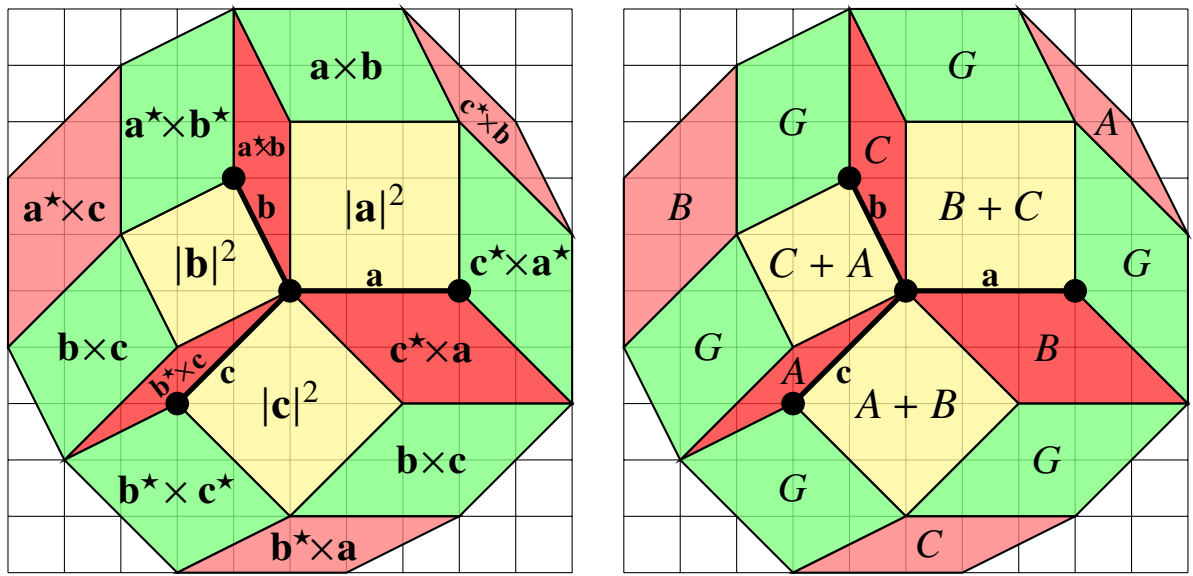

~\\
{\bf  Proof:}    Let us start with three spinors between three tangent disks.
Build the network of spinor as shown in Figure \ref{fig:distribution}.
The signs of the clockwise oriented spinors $\mathbf a$, $\mathbf b$, and $\mathbf c$ are chosen so that 
$$
\mathbf a + \mathbf b + \mathbf c = 0
$$
Figure~\ref{fig:labeled} shows the tessellation based on these spinors. 
\\\\
{\bf 1.}	The red parallelograms are the curvatures of the Descartes by the virtue of \eqref{eq:thm4}:
$$
                A = \mathbf b\times \mathbf c^\star ,  \qquad   B 
                     = \mathbf c \times \mathbf a^\star ,  \qquad C = \mathbf a \times \mathbf b^\star 
$$
The areas of yellow squares that represent the squares of the the spinors:  
$$
                   |\mathbf a|^2 = B+C, \qquad  |\mathbf b|^2 = C+A, \qquad  |\mathbf c|^2 = A+B\,,
$$
This is consistent with  
$$
       |\mathbf a|^2 = \mathbf a\times \mathbf a^\star 
       = (-\mathbf b-c)\bigcdot \mathbf a^\star = -\mathbf b\times \mathbf a^\star - \mathbf c\times \mathbf a^\star 
        = \mathbf a^\star \times \mathbf b + \mathbf c^\star \times \mathbf a = B+C.
$$    
Using \eqref{eq:thm4} again, one finds the curvature of the fourth inscribed circle:
$$
D \quad = \quad   (\mathbf a+\mathbf c^\star)\times ( \mathbf b+\mathbf a^\star)  \quad = \quad 
                      \underbrace { ~\mathbf a\times \mathbf b\phantom{\Big|}}_{\hbox{\smalll\sf green}\atop G}  \ + \
                      \underbrace{~\mathbf a\times \mathbf a^\star\phantom{\Big|}}_{\hbox{\smalll\sf square} \atop B+C} \ + \
                      \underbrace{~\mathbf c^\star\times \mathbf b\phantom{\Big|}}_{\hbox{\smalll\sf thin red} \atop A } \ + \ 
                    \underbrace{~\mathbf c^\star\times \mathbf a^\star\phantom{\Big|}}_{\hbox{\smalll\sf green} \atop G }
$$   
which indeed may be composed  into a butterfly shape as in Figure~\ref {fig:3butterflies}.

\noindent
{\bf 2.}	All green tiles have the same area:
$$
               \mathbf a\times \mathbf b  =  \mathbf b\times \mathbf c  =  \mathbf  c\times \mathbf a  
               =  \hbox{curvature of the mid-circle} (A,B,C)
$$
which follows from
$$
\mathbf a\times \mathbf b = \mathbf a\times (-\mathbf c-\mathbf a) 
                                      = -\mathbf a\times \mathbf c - \mathbf a\times \mathbf a 
                                      = \mathbf c\times \mathbf a
$$ 
and from $\mathbf a\times \mathbf b = \mathbf a^\star\times \mathbf b^\star$.
Using \eqref{eq:thm42} for two spinors into $D$ we obtain the curvature of the mid-circle $(ABC)$:
$$
               \mathbf b^\star\bigcdot \mathbf c \ = \ \mathbf b \times \mathbf c 
               %  =   c\times a  =  \hbox{curvature of the mid-circle} (A,B,C)
$$
which is one of the green tiles (see Figure \ref{fig:labeled}).
As to the other three other mid-circles, one of them is
$$
(ABD) = (\mathbf c+\mathbf b^\star)\bigcdot \mathbf c 
           = \mathbf c\bigcdot \mathbf c + \mathbf b^\star\bigcdot \mathbf c 
            = |\mathbf c|^2 + \mathbf b \times \mathbf c
$$
which may be recognized as a yellow tile plus green.
\\[-5pt]
~

\noindent
{\bf 3.}  As to the other solution for the Descartes problem, 
we use again \eqref{eq:thm4} for two spinors into $D'$:
$$
\begin{array}{rl}
(\mathbf a+\mathbf b^\star) \times (\mathbf c+\mathbf a^\star)
\ \ &= \ \  \mathbf a\!\times\! \mathbf c   +  \mathbf  a\!\times\! \mathbf a^\star
               \ +\  \mathbf b^\star\!\!\times\! \mathbf c \   + \ \mathbf b^\star \!\!\times \! \mathbf a^\star\\[5pt]
\ \ &= \  \  |\mathbf a|^2 \ + \ \mathbf b^\star\times \mathbf c \ - \  \mathbf a\times \mathbf b - \mathbf c\times \mathbf a
\end{array}
$$
which is a yellow tile plus  a red tile minus two green tiles.
The symmetry mid-circles for the Descartes configuration $(ABCD')$ can be obtained similarly.
For instance
$$
\begin{array}{rlll}
(BCD') \  = \ - (\mathbf a+\mathbf b^\star)	\bigcdot(\mathbf c+\mathbf a^\star) 
      &=& - \mathbf a\bigcdot \mathbf c - \mathbf b^\star\bigcdot \mathbf c - \mathbf b^\star\bigcdot \mathbf a^\star 
          &\qquad(\hbox{since }\mathbf a\bigcdot \mathbf a^\star = 0)\\
      &=&  \mathbf a\bigcdot(\mathbf a+\mathbf b) -  \mathbf b\times \mathbf c -\mathbf a\bigcdot \mathbf b \\
      &=& |\mathbf a|^2 - \mathbf b\times \mathbf c
\end{array}
$$
which a difference between a yellow and a green tile.
\QED
\\

%$$
%\begin{array}{rlll}
%(a+b^\star)\cdot(b+c^\star)	 &=& a\cdot b + a\cdot c^\star + b^\star\cdot c^\star \\
%                                 &=& -(b+c)\cdot b +c\times a +b\cdot c \\
%                                 &=& -|b|^2 +c\times a
%                                 \end{array}
%$$
%
%$$
%\begin{array}{rlll}
%(b+c^\star)\cdot(c+a^\star)	 &=& b\cdot c + b\cdot a^\star + c^\star\cdot a^\star \\
%                                 &=& -(c+a)\cdot c +a\times b +c\cdot a \\
%                                 &=& -|c|^2 +a\times b
%                                 \end{array}
%$$
%
%$$
%\begin{array}{rlll}
%(c+a^\star)\cdot(a+b^\star)	 &=& c\cdot a + c\cdot b^\star + a^\star\cdot b^\star \\
%                                 &=& -(a+b)\cdot a +b\times c +a\cdot b \\
%                                 &=& -|a|^2 +b\times c
%                                 \end{array}
%$$
%

An obvious corollary to this is that if spinors are integral, so is the Descartes configuration.
%\\\\

%------------------------------------------------------------------
\section{Coda: Integral Descartes configurations from spinors}

Stripping the content of the above construction from its geometric interpretation,  we get effectively a method 
of obtaining integral Descartes configurations (and consequently integral Apollonian disk packings) 
from four arbitrary integers.
The core is the curl property \eqref{eq:curl}: 
%$$
\begin{equation}
\label{eq:again}
                                \mathbf a+\mathbf b+\mathbf c = \mathbf 0
\end{equation}
%$$
The implied curvatures of the disks in the Descartes configuration, $A$, $B$,  $C$, and the two alternative curvatures of the fourth, $D_1$ and $D_2$,
are:
%$$
\begin{equation}
\label{eq:short}
\begin{array}{rl}
A  &= \ \ -\mathbf b\bigcdot \mathbf c\\
B  &=\ \ -\mathbf c\bigcdot \mathbf a \\
C &=\ \ -\mathbf a\bigcdot \mathbf b\\
D_1+D_2 &=\ \ |\mathbf a|^2  +|\mathbf b|^2 + |\mathbf c|^2 \\
D_1-D_2 &=\ \ 4\mathbf a\times \mathbf b = 4 \mathbf  b\times \mathbf c = 4\mathbf c\times \mathbf a
\end{array}
\end{equation}
%$$
Reducing the above to two spinors (4 integers) with \eqref{eq:again} leads to the claim.
\\[5pt]
{\bf Proposition 7:} Let $\mathbf a$ and $\mathbf b$ be arbitrary integral vectors in $\mathbb Z^2$.  Then the following 
integers satisfy the Descartes Diophantine equation \eqref{eq:Descartes}:
 %$$
\begin{equation}
\label{eq:reduced}
\begin{array}{rl}
A  &= \ \  |\mathbf b|^2 +\mathbf a\bigcdot \mathbf b\\
B  &=\  \   |\mathbf a|^2 +\mathbf a\bigcdot \mathbf b \\
C  &=\ \ -\mathbf a\bigcdot \mathbf b\\
D  &=  \ \ |\mathbf a|^2+|\mathbf b|^2  \; +\; \mathbf a\bigcdot \mathbf b 
%+
\; \pm \; 2\, 
\mathbf a\times \mathbf b
\end{array}
\end{equation}
%$$ 

~

A different parametrization of Descartes configurations from the spinors 
based on the divergence property \eqref{eq:thm5a} will be presented in a separate paper.

~

%---------------------------------------------
%\section*{Note on software}
\noindent
{\bf Note on software:}
Most of the figures were made with Tikz \cite{tikz}
except Figure 6--8, which were made with Cinderella \cite{cinderella}. 
%Figures \ref{fig:MainExample} and \ref{fig:more} were made with the help of the software Cinderella \cite{cinderella}.
%The remaining figures are made with TikZ \cite{tikz}. 

%Bibliography-------------------------------------------------------------------------------

\end{document}